\documentclass[leqno,12pt]{article}
\usepackage{latexsym}
 \usepackage{graphicx}
\usepackage{amsmath}
\usepackage{amssymb}
\usepackage{mathtools}
\usepackage{cite}
\usepackage{color}
\usepackage{subfig}
\usepackage[margin=1.3in, top=1.3in, bottom=1.3in]{geometry}

\usepackage[colorlinks,linkcolor=black,citecolor=blue,urlcolor=blue,bookmarks=false,hypertexnames=true]{hyperref} 

\newcommand{\nc}{\newcommand}
\nc{\nt}{\newtheorem}
\nt{thm}{Theorem}[section]
\nt{cor}[thm]{Corollary}
\nt{prop}[thm]{Proposition}
\nt{lem}[thm]{Lemma}
\nt{defn}[thm]{Definition}
\nt{rem}[thm]{Remark}
\nt{exa}[thm]{Example}
\nt{ass}[thm]{Assumption}
\nt{alg}[thm]{Algorithm}
\nt{conj}[thm]{Conjecture}
\nc{\ip}[2]{\mbox{$\langle #1,#2 \rangle$}}
\nc{\pf}{\noindent{\bf Proof\ \ }}
\nc{\finpf}{\hfill{$\Box$}\linespace}
\nc{\linespace}{\vspace{\baselineskip} \noindent}
\nc{\R}{{\mathbf R}}
\nc{\oR}{\overline{\R}}
\nc{\M}{\mathcal M}

\nc{\Rn}{{\mathbf R}^n}
\nc{\inT}{\mbox{\rm int}\,}
\nc{\cl}{\mbox{\rm cl}\,}

\def\tto{\;{\lower 1pt \hbox{$\rightarrow$}}\kern -12pt
           \hbox{\raise 2.8pt \hbox{$\rightarrow$}}\;}
\newenvironment{myequation}{\setcounter{equation}{\value{thm}}
   \begin{equation}}{\addtocounter{thm}{1}\end{equation}}

\nc{\bmye}{\begin{myequation}}
\nc{\emye}{\end{myequation}}

\nc{\red}{\textcolor{red}}
\nc{\blue}{\textcolor{blue}}
\usepackage[normalem]{ulem} 

\begin{document}
\title{
Identifiability, the KL property in metric spaces, and subgradient curves
}
\author{
\and
A.S. Lewis
\thanks{ORIE, Cornell University, Ithaca, NY.
\texttt{people.orie.cornell.edu/aslewis} 
\hspace{2cm} \mbox{~}
Research supported in part by National Science Foundation Grant DMS-2006990.}
\and
Tonghua Tian
\thanks{ORIE, Cornell University, Ithaca, NY.  
\texttt{tt543@cornell.edu}}
}
\date{\today \\ \bigskip
{\normalsize Communicated by Michael Overton}
}
\maketitle

\begin{abstract}
Identifiability, and the closely related idea of partial smoothness, unify classical active set methods and more general notions of solution structure.  Diverse optimization algorithms generate iterates in discrete time that are eventually confined to identifiable sets.  We present two fresh perspectives on identifiability.  The first distills the notion to a simple metric property, applicable not just in Euclidean settings but to optimization over manifolds and beyond; the second reveals analogous continuous-time behavior for subgradient descent curves.  The Kurdyka-{\L}ojasiewicz property typically governs convergence in both discrete and continuous time:  we  explore its interplay with identifiability.
\end{abstract}
\medskip

\medskip


\section{Introduction}
Contemporary optimization involves notions of structure, such as sparsity and rank, that have prompted a re-examination of classical active-set philosophy.  Early generalizations and terminology such as \cite{Flam92,Burke-More88,Burke90,Wright} motivated the idea of an {\em identifiable set} for an objective function $f$ over Euclidean spaces \cite[Definition 3.10]{ident}.

Consider, for example, the nonsmooth nonconvex function on $\R^2$ defined by
\bmye \label{simple1}
f(x) ~=~ 5 |x_2-x_1^2| + x_1^2.
\emye
The set $\M = \{ x : x_2=x_1^2 \}$
is identifiable at the minimizer $\bar x = 0$ because, in the language of subgradients \cite{VA}, a quick calculation shows
\bmye \label{identifiable-example}
\mbox{dist}\big(0,\partial f(x^k)\big) \to 0 \qquad \Rightarrow \qquad x^k \in \M~
\mbox{eventually}.
\emye

A wide variety of algorithms for nonsmooth or constrained optimization generate convergent sequences of iterates that behave well with respect to some merit function.  Specifically, some corresponding sequence of subgradients of the merit function converges to zero. Consequently, the algorithm identifies some associated structure --- in example (\ref{simple1}), the set $\M$.  Examples  include classical projected gradient or subgradient methods \cite{Burke-More88,Wright,Flam92}, the proximal point method \cite{ident_active} and composite generalizations such as the proximal gradient method \cite{prx_lin}, alternating projection and Gauss-Seidel schemes \cite{desc_semi,bolte_complex_kur}, some bundle methods \cite{noll-kl}, and the broad majorization-minimization framework \cite{bolte-pauwels-mm}, including sequential quadratic programming.  Interesting contemporary examples include \cite{fb_ps,activ_ident}. 

Within optimization and its interface with machine learning, interest in non-Euclidean settings has grown rapidly, supported by expositions like \cite{absil,boumal2022intromanifolds}.  In particular, \cite{zhang-sra} proved complexity bounds for the projected subgradient approach to geodesically convex optimization on Hadamard manifolds.  Earlier work had already extended the proximal point philosophy far beyond Euclidean settings, in particular to Hadamard manifolds: \cite{ferreira-oliveira-02,li-lopez-martin-marquez-09,bento-ferreira-oliveira-10} study geodesically convex and nonconvex optimization, and alternating projection schemes.  Extensions to more general geodesic spaces include \cite{jost,bacak-searston-sims-12,bacak} (in the case of nonpositive curvature) and \cite{lewis-lopez-nicolae-22}.  A reappraisal of the fundamentals of identifiability from a purely metric perspective is therefore appealing. 

\begin{figure}
    \centering
    \includegraphics[width=0.6\textwidth]{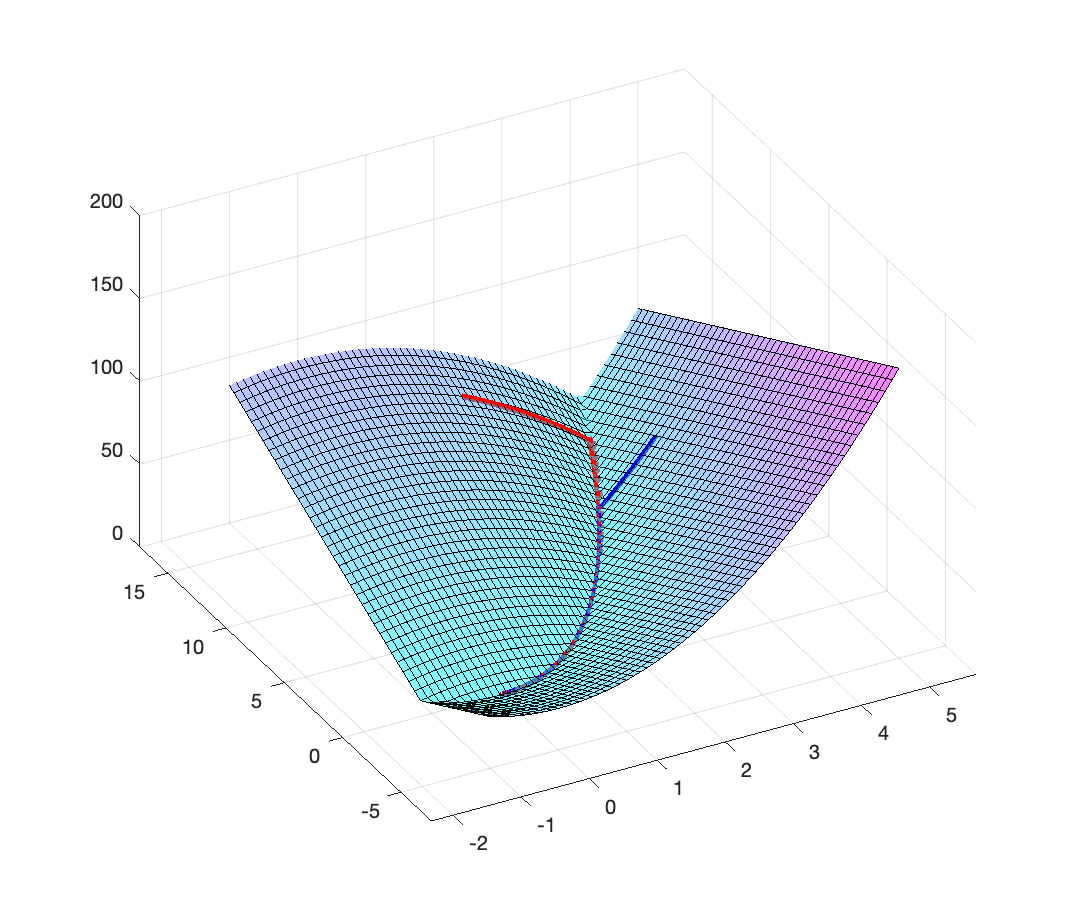}
    \caption{Subgradient curves for the objective $5 |x_2-x_1^2| + x_1^2$.}
    \label{fig:subcurve}
\end{figure}

A complementary perspective on identifiability arises through analogues of iterative algorithms like the proximal point method but in continuous time \cite{brezis}.  The same identification behavior illustrated above in discrete time also manifests itself in continuous time.  In the example, as we shall see later, locally absolutely continuous curves $x \colon \R_+ \to \Rn$ with initial points $x(0)$ near zero and satisfying
\[
x'(t) \in -\partial f\big(x(t)\big) \qquad \mbox{for almost all times}~ t>0
\]
(following standard variational-analytic terminology \cite{VA})
always converge to zero and are eventually confined to $\M$.  Figure \ref{fig:subcurve} shows plots of two randomly initialized examples of these {\em subgradient curves}. 

We pursue three main themes around the idea of identifiability for an objective $f$.  First, we reframe the idea in purely {\em metric} terms, using the {\em slope} $|\nabla f|$, which measures instantaneous rate of decrease \cite{giorgi}.  For example, property (\ref{identifiable-example}) becomes
\[
|\nabla f|(x^k) \to 0 \qquad \Rightarrow \qquad x^k \in \M~
\mbox{eventually}.
\]
By switching from the sophisticated dual perspective of subgradient vectors to the simple primal scalar perspective of the slope, we highlight the intuitive nature of identifiability and extend its applicability to more general settings like optimization over manifolds \cite{absil,boumal2022intromanifolds}.  This setting also allows us to highlight a new linear growth tool (Theorem \ref{linear}), on which our whole development fundamentally depends. 
Secondly, in the Euclidean setting, we demonstrate quite generally that subgradient curves are eventually confined to identifiable manifolds.  In the Euclidean setting, our metric definition of identifiability is equivalent to the usual subgradient-based  definition, and closely related to {\em partial smoothness} \cite{Lewis-active}.  Thirdly, motivated by the crucial role of the ``Kurdyka-\L ojasiewicz property'' in the convergence of optimization algorithms and dynamics \cite{loja,Lewis-Clarke,desc_semi,bolte_complex_kur}, including those in metric-space settings \cite{hauer-mazon}, we show how the KL property is controlled entirely by the KL property for the restriction of the objective to any identifiable set.  Optimization problems we encounter in practice are semi-algebraic or, more generally, tame \cite{tame_opt}, and for generic problems instances in such settings, local minimizers lie on identifiable analytic manifolds on which the objective function is analytic \cite{gen,gen_exist}.  The KL property for analytic functions on analytic manifolds is classical \cite[Section 9]{kurdyka-mostowski-parusinski}, so our results  apply, shedding a fresh light on the KL property in nonsmooth optimization.   

We present our development in two parts.  Part \ref{part:1} --- Sections \ref{sec:Identifiability} through \ref{KL} --- studies identifiability and its interplay with the KL property in metric spaces;  Part~\ref{part:2} --- Sections \ref{four} through \ref{sec:pln} --- focuses on the special case of Euclidean space.  

Part \ref{part:1} begins in Section \ref{sec:Identifiability} by reviewing the notion of slope in metric spaces, and its relationship with subgradients in the Euclidean case.  We introduce the idea of a critical sequence, noting the prevalence of such sequences in optimization as proximal images of minimizing sequences.  Thus motivated, Section \ref{sec:identifiable} presents our new metric definition of identifiability:  sets in which all critical sequences are eventually contained.  In particular, expanding on a Euclidean precedent \cite{ident}, we study minimal identifiable sets.  In Section \ref{sec:linear}, in striking contrast to previous subgradient-based developments, we deploy the Ekeland variational principle to explain the fundamental geometry of linear objective growth around identifiable sets in complete metric spaces. Our new linear growth result explains, in broad generality, how objective growth around minimizers is governed entirely by growth on the identifiable set.  As Section \ref{KL} explains, a particularly important growth condition is the Kurdyka-\L ojasiewicz property, which we newly reframe in terms of identifiability.  We prove that the KL property is inherited from the analogous property on an identifiable set, a useful tool for understanding the property in practice.

Moving to the Euclidean setting of Part \ref{part:2}, we begin in Section \ref{four} with a quick review of some variational analysis.  Section \ref{max} develops the basic example of ``max functions'':  pointwise maxima of finitely many continuously differentiable functions.   Section~\ref{sec:ps} studies identifiability from a variational-analytic perspective, throwing a fresh and intuitive light on the equivalence between identifiable manifolds and partial smoothness \cite{ident}. We return in particular to the KL property, now relying on analytic settings.
In Section \ref{sec:Subgradient} we move to the continuous time setting, presenting a result that explains for the first time the identification behavior illustrated in Figure~\ref{fig:subcurve}. We end in Section \ref{sec:pln} with a proof of this result.

\part{Identifiability in metric spaces}\label{part:1}

\section{Critical sequences in metric spaces} \label{sec:Identifiability}
Our setting is a metric space $(X,d)$.  For any point $x \in X$, the function $d_x \colon X \to \R$ is defined by $d_x(y) = d(x,y)$ for points $y \in X$.  On $X$, we consider a function $f$ which, unless otherwise stated, takes values in $(-\infty,+\infty]$.  We focus on a point $\bar x$ in the domain 
\[
\mbox{dom}\,f ~=~ \{x : f(x)<+\infty\}.
\] 
We are mainly concerned with points $\bar x$ that are local minimizers for $f$.  First, however, we recall a fundamental definition \cite{giorgi} leading to a weaker property.

\begin{defn}
{\rm
For a function $f$ on a metric space, the {\em slope} $|\nabla f|(\bar x)$ at a point 
$\bar x \in \mbox{dom}\,f$ is zero if 
$\bar x$ is a local minimizer, and otherwise is the quantity
\[
\limsup_{\bar x \ne x \to \bar x} \frac{f(\bar x) - f(x)}{d(\bar x, x)}.
\] 
At points outside of the domain, the slope is $+\infty$.
}
\end{defn}
The notation is suggestive of the easy case of a differentiable function on a Euclidean space, when 
$|\nabla f|(x) = |\nabla f(x)|$ \cite[Observation 4.3]{giorgi}.  Notice that the slope may be infinite at points in the domain:  for example, $|\nabla \sqrt[3]{\cdot}|(0) = + \infty$.  An easy exercise shows an equivalent definition:
\[
|\nabla f|(\bar x) ~=~ \inf \{\delta > 0 : \mbox{$\bar x$ is a local minimizer for $f + \delta d_{\bar x}$} \}.
\]
(The infimum of the empty set is defined to be $+ \infty$.) By definition, the slope is zero at any local minimizer.  However, the converse may fail:  for example, on $\R$, the function $f(x)=-x^2$ has slope zero at zero but is not locally minimized there.  The following result \cite[Proposition 3.3]{ioffe-variational} collects some simple tools for the slope.

\begin{prop} \label{slope-tools}
Functions $f$ and $g$ on a metric space always satisfy 
\[
|\nabla(f+g)|(x) ~\le~ |\nabla f|(x) + |\nabla g|(x) \qquad \mbox{and} \qquad 
|\nabla(\delta f)|(x) ~=~ \delta |\nabla f|(x)
\]
for all points $x$ and scalars $\delta > 0$.  Furthermore, if $f$ is $L$-lipschitz, then 
$|\nabla f|(x) \le L$.
\end{prop}

\noindent
By Proposition \ref{slope-tools}, all points $x,y$ in the metric space $X$ satisfy $|\nabla d_x|(y) \le 1$.  When $X$ is a {\em length space}, meaning that the distance $d(x,y)$ is always the infimum of the lengths of curves joining $x$ and $y$, we can say more.

\begin{thm}[Slope characterization of length spaces \cite{aze-corvellec}] \label{length}
A complete metric space $(X,d)$ is a length space if and only if all distinct points $x,y \in X$ satisfy 
\mbox{$|\nabla d_x|(y) = 1$}.
\end{thm}

\noindent
Any connected Riemannian manifold with its Riemannian distance is a length space \cite[Theorem 2.55]{lee-riemannian}.  We also observe the following easy inequality.

\begin{prop} \label{square}
On any metric space, the squared distance function satisfies 
\[
|\nabla d_x^2|(y) ~\le~ 2d(x,y) \qquad \text{and} \qquad |\nabla (-d_x^2)|(y) ~\le~ 2d(x,y)
\]
for all points $x,y \in X$.
\end{prop}

\pf
Since
\begin{eqnarray*}
\limsup_{y \ne z \to y} \frac{d_x^2(y) - d_x^2(z)}{d(y,z)} 
&=& \limsup_{y \ne z \to y}  \frac{\big(d(y,x) - d(z,x)\big)\big(d(y,x) + d(z,x)\big)}{d(y,z)} \\
& \le & \limsup_{y \ne z \to y} \big(d(y,x) + d(z,x)\big) ~=~ 2d(y,x),
\end{eqnarray*}
the first inequality follows. The second inequality follows similarly.
\finpf

A property even weaker than the slope being zero is the criticality idea that we define next.  Following \cite{VA}, for any function $f$ and any point $\bar x \in X$, we write $x_r \to_f \bar x$ to denote {\em $f$-attentive convergence\/}:  $x_r \to \bar x$ and 
$f(x_r) \to f(\bar x)$.  Furthermore, for any function $g$, we define the quantity
\[
\liminf_{x \to_f \, \bar x} g(x) ~=~ 
\inf \{ \liminf_r g(x_r) : x_r \to_f \bar x\},
\]
which by definition is no larger than $g(\bar x)$.

\begin{defn} \label{def:critical}
{\rm
Given a function $f$, a {\em critical sequence} for a point $\bar x \in \mbox{dom}\,f$ is a sequence 
$x_r \to_f \bar x$ such that $|\nabla f|(x_r) \to 0$, and $\bar x$ is {\em critical} if such a sequence exists.
}
\end{defn}
For example, any sequence $x_r \to x$ in $X$ is critical for the function $d_x^2$, by Proposition~\ref{square}.  On $\R$, the function $f(x) = \min\{x,0\}$ has slope $1$ at the point $0$, but that point is critical because the local minimizers $\frac{1}{r}$ for $r=1,2,3,\ldots$ form a critical sequence.  The notion of a critical point is related to the following more robust version of the slope \cite[Definition 8.4]{ioffe-variational}.

\begin{defn}
{\rm
The {\em limiting slope\/} of $f$ at a point $\bar x \in \mbox{dom}\,f$ is the quantity
\[
\overline{|\nabla f|}(\bar x) ~=~ \liminf_{x \to_f \, \bar x} |\nabla f|(x).
\]
}
\end{defn}
With this terminology, critical points are exactly those where the limiting slope is zero.  Notice that the limiting slope is never larger than the slope.

Our development of identifiability in metric spaces studiously avoids the use of subgradients.  However, in the special case when the underlying space is Euclidean, readers familiar with standard variational analysis (for example \cite{VA}) will find the following formula \cite[Proposition 8.5]{ioffe-variational} illuminating.  It relates the limiting slope of a closed function with the usual subdifferential $\partial f$ in the terminology of \cite{VA}.

\begin{prop}[Distance formula]~ \label{distance-formula}
On $\Rn$, any closed function $f$ and any point $x \in \mbox{\rm dom}\, f$ satisfy
$\overline{|\nabla f|}(x) = \mbox{\rm dist}\big( 0 , \partial f(x) \big)$.
\end{prop}

\noindent
In Euclidean space, critical points are therefore those at which zero is a subgradient.

A wide variety of optimization procedures generate critical sequences, because the procedures involve projection or proximal operations corresponding to closed functions $f \colon X \to (-\infty,+\infty]$:
\bmye \label{proximal-operator}
\mbox{prox}_f (x) ~=~ \mbox{argmin} \Big\{f + \frac{1}{2} d_x^2 \Big\}.
\emye
In Euclidean settings, this observation is very familiar \cite{Burke-More88,Wright,Flam92,ident_active,prx_lin,desc_semi,bolte_complex_kur,noll-kl,fb_ps,activ_ident}.  However, the phenomenon persists in non-Euclidean settings such as Hadamard manifolds 
\cite{zhang-sra,ferreira-oliveira-02,li-lopez-martin-marquez-09,bento-ferreira-oliveira-10} and more general geodesic spaces \cite{jost,bacak-searston-sims-12,bacak,lewis-lopez-nicolae-22}.  We next discuss this property from a purely metric perspective.  

The proximal operator  is single-valued for closed convex functions $f$ on a Hilbert space (by Minty's Theorem \cite{minty}) and in some more general metric spaces \cite{jost,bacak}.  However, for general spaces and closed functions, the value $\inf\{f + \frac{1}{2}d_x^2\}$ can be $-\infty$, and even when finite can be unattained.  We do, however, have the following routine result in {\em proper} metric spaces (meaning that all closed balls are compact).

\begin{lem}
On a proper metric space $X$, if a closed function $f \colon X \to (-\infty,+\infty]$ is not identically $+\infty$ and is bounded below, then the operator $\mbox{\rm prox}_f$ is everywhere nonempty-valued.
\end{lem}

\pf
Consider any lower bound $\alpha$ for $f$ and any point $z \in X$ for which $f(z)$ is finite.  Given any point $x \in X$, the function $g = f + \frac{1}{2}d_x^2$ is closed, so the level set 
$L = \{y : g(y) \le g(z)\}$ is closed.  Notice that $L$ contains $z$, and furthermore is bounded, since all points $y \in L$ satisfy
$\frac{1}{2}d_x^2(y) \le g(z) - \alpha$.  Hence $L$ is nonempty and compact, so the set $\mbox{prox}_f (x) = \mbox{argmin}_L g$ is nonempty.
\finpf

The proximal operator has the obvious property that any minimizer of the function $f$ is a {\em fixed point\/}:  $x \in \mbox{prox}_f(x)$, and more generally,  $x \in \mbox{prox}_{\frac{1}{\alpha}f}(x)$ for any constant $\alpha > 0$.  In turn, $f$ has slope zero at any such fixed point.  In Euclidean space this follows from the property that $y \in \mbox{\rm prox}_f(x)$ implies $x-y \in \partial f(y)$.  More generally, it is a consequence of the following easy tool.

\begin{lem} \label{lem:proximal}
For any closed function $f$ on a metric space, any points $x \in \mbox{\rm dom}\, f$ and $y \in \mbox{\rm prox}_f(x)$ satisfy 
$|\nabla f|(y) \le d(x,y)$.  Consequently, $f$ has slope zero at any fixed point of $\mbox{\rm prox}_f$.
\end{lem}

\pf
Propositions \ref{slope-tools} and \ref{square} and the fact that $\big|\nabla \big(f + \frac{1}{2} d_x^2\big)\big|(y) = 0$ imply
\[
|\nabla f|(y) ~\le~ \big|\nabla \big(f + \frac{1}{2} d_x^2\big)\big|(y) + \big|\nabla \big(- \frac{1}{2} d_x^2\big)\big|(y) ~\le~ d(x, y).
\]
The final claim follows by setting $y=x$.
\finpf

We next show that the proximal operator, when nonempty-valued, allows us to associate a critical sequence with any minimizing sequence.  Special cases of this mechanism, which underlies a great deal of  the identification literature that we discussed after equation~(\ref{proximal-operator}), date back at least to \cite{bertsekas}, although this simple metric version seems a new perspective.  

\begin{prop}[Proximal images of minimizing sequences]\label{prop:proximal}
For a closed function $f$ on a metric space and constants $\alpha > \beta > 0$, consider any point 
$\bar x$ that minimizes $f$, or more generally, is a fixed point of the proximal operator 
$\mbox{\rm prox}_{\frac{1}{\beta} f}$.  Then,  for any sequence of points $x_r \to_f \bar x$, any sequence of points $y_r \in \mbox{\rm prox}_{\frac{1}{\alpha} f} (x_r)$ for $r=1,2,3,\ldots$ is critical for $\bar x$.
\end{prop}

\pf
We have
\begin{eqnarray*}
f(x_r) - f(\bar x)
& \ge &
f(y_r) + \frac{\alpha}{2}  d(y_r,x_r)^2 - f(\bar x) 
~ \ge ~
\frac{\alpha}{2} d(y_r,x_r)^2 - \frac{\beta}{2} d(y_r,\bar x)^2  \\
& \ge &
\frac{\alpha}{2} d(y_r,x_r)^2 - \frac{\beta}{2} \big( d(y_r,x_r) + d(x_r,\bar x) \big)^2 \\
& \ge &
\frac{\alpha-\beta}{2} d(y_r,x_r)^2 - \beta d(y_r,x_r)d(x_r,\bar x) - \frac{\beta}{2}  d(x_r,\bar x)^2.
\end{eqnarray*}
Rearranging gives
\[
f(x_r) - f(\bar x) + \frac{\alpha\beta}{ 2(\alpha-\beta)} d(x_r,\bar x)^2
~\ge~
\frac{\alpha-\beta}{2} \Big( d(y_r,x_r) - \frac{\beta}{\alpha-\beta} d(x_r,\bar x) \Big)^2.
\]
Since the left-hand side converges to zero, so does the right-hand side, so we deduce $d(y_r,x_r) \to  0$, and hence $y_r \to \bar x$.  Since the function $f$ is closed, we know 
$\liminf f(y_r) \ge f(\bar x)$.  But since
$f(y_r) + \frac{\alpha}{2} d(y_r,x_r)^2 \le f(x_r)$
for each $r$, we also know  $\limsup f(y_r) \le f(\bar x)$, so $f(y_r ) \to f(\bar x)$.  Finally, Proposition \ref{slope-tools} and Lemma~\ref{lem:proximal} complete the proof:
$|\nabla f|(y_r) \le \alpha d(y_r,x_r) \to 0$.
\finpf

Although minimizers are always fixed points of the proximal operator, points where the slope is zero may not be.  For example, the real function $f(x) = -|x|^{\frac{3}{2}}$ has slope zero at the point $\bar x=0$, but $\bar x$ is not a fixed point of $\mbox{\rm prox}_{\frac{1}{\beta} f}$ for any value $\beta > 0$.  However, when the space $X$ is Euclidean, for well-behaved functions $f$, any critical point is a fixed point of $\mbox{\rm prox}_{\frac{1}{\beta} f}$ for all sufficiently large $\beta > 0$.  Specifically, following the standard terminology of \cite{VA}, this property holds if $f$ is both ``prox-regular'' and ``prox-bounded'' \cite[Proposition 13.37]{VA}.  

\section{Identifiable sets in metric spaces} \label{sec:identifiable}
In the previous section we saw how minimizing sequences are often associated with critical sequences.
As we discussed in the introduction, many optimization problems possess inherent structure around critical points that restricts the possible corresponding critical sequences.  To formalize this phenomenon, and unify many precedents in the literature (originating with \cite{Calamai-More87}), \cite{ident} used the language of subgradients in Euclidean spaces to develop the notion of an ``identifiable set'', and to study when such sets are ``locally minimal'' \cite[Definitions 3.10 and 3.11]{ident}.  The terminology originated with \cite{Burke-More88}.  Our contribution is a fresh view and extension of this theory from a simple metric space perspective.

\begin{defn} \label{def:identifiability}
{\rm
Given a function $f$ on a metric space $X$, consider a point $\bar x \in \mbox{dom}\, f$ and a set $\M \subset X$ containing $\bar x$.  
\begin{itemize}
\item
$\M$ is {\em identifiable} at $\bar x$ if its complement contains no critical sequence for $\bar x$. 
\item
$\M$ is {\em strongly critical} at $\bar x$ if all sequences $(x_r)$ in $\M$ satisfying 
$x_r \to_f \bar x$ are critical for $\bar x$.
\item
The {\em modulus of identifiability} for $\M$ at $\bar x$ is
\bmye \label{rate}
\liminf_{x \to_f \, \bar x,~ x \not\in \M} |\nabla f |(x).
\emye
\end{itemize}
}
\end{defn}

\noindent
Notice that $\M$ is identifiable if and only if the modulus of identifiability is strictly positive.  Furthermore, if the point $\bar x \in \mbox{dom}\, f$ is not critical, then every set containing 
$\bar x$ is identifiable, so our interest in identifiable sets focuses on critical points.  Identifiable sets always exist.  In particular, for any radius $\epsilon > 0$, setting 
$B_{\epsilon}(\bar x) = \{x \in  X : d(x,\bar x) \le \epsilon \}$, the {\em $f$-attentive neighborhood}
\bmye \label{attentive}
B_{\epsilon}^f(\bar x) ~=~ \{x \in  B_\epsilon(\bar x) : |f(x) - f(\bar x)| \le \epsilon\}
\emye 
is identifiable. The following tools are easy to verify.

\begin{prop} \label{finite}
On a metric space, at any point in the domain of a function, finite intersections of identifiable sets are identifiable, and finite unions of strongly critical sets are strongly critical.
\end{prop}

\begin{prop} \label{strongly}
On any metric space, consider functions $f \ge g$, a point \mbox{$\bar x \in \mbox{\rm dom}\, f$} and a set $\M$ containing $\bar x$.  Suppose that $f$ and $g$ are equal on $\M$ around $\bar x$.  At $\bar x$, if $\M$ is strongly critical for $g$, then $\M$ is also strongly critical for $f$.
\end{prop}

\pf
For any sequence $x_r \to_f \bar x$ in $\M$  since $g \le f$ and $g(x_r) = f(x_r)$, strong criticality for $g$ implies $0 \le |\nabla f|(x_r) \le |\nabla g|(x_r) \to 0$, so $(x_r)$ is critical for $f$.
\finpf

\begin{rem} \label{rem:locality}
{\rm
If a set $\M$ is identifiable at a point $\bar x$, then it is easily seen to be identifiable at every point in the intersection of $\M$ with some $f$-attentive neighborhood of $\bar x$. Moreover, the modulus of identifiability $\sigma_x$, as a function of $x \in \M$, is lower semicontinuous with respect to $f$-attentive convergence:
\[
\liminf_{x \to_f \, \bar x, ~ x \in \M} \sigma_x ~\ge~ \sigma_{\bar x}.
\]
}
\end{rem}

Identifiable sets are not unique:  any superset of an identifiable set at a point $\bar x$ is also identifiable.  Identifiability is a local property ($f$-attentively):  if a set is identifiable, so is its intersection with any $f$-attentive neighborhood.  Our interest in strongly critical sets stems from the fact that smaller identifiable sets are more informative tools than larger ones.   To be precise, any strongly critical identifiable set at $\bar x$ must be locally minimal (in an $f$-attentive sense) among identifiable 
sets. A subgradient-based version of the following result appeared in \cite[Proposition 3.12]{ident}.

\begin{prop}[Locally minimal identifiable sets] \label{minimality}
Given a function $f$ on a metric space, suppose that a set $\M$ is identifiable at a point $\bar x$.  Then $\M$ is strongly critical at $\bar x$ if and only if the following property holds:
\bmye \label{minimality2}
\left\{
\begin{array}{l}
\mbox{For any identifiable set $\M'$ at $\bar x$, there exists an $f$-attentive} \\ 
\mbox{neighborhood $B_{\epsilon}^f(\bar x)$ with $\epsilon > 0$ satisfying 
$\M \cap B_{\epsilon}^f(\bar x) \subset \M'$.}
\end{array}
\right.
\emye
\end{prop}

\pf
Suppose that the set $\M$ is strongly critical at $\bar x$ but that property (\ref{minimality2}) fails for some set $\M'$ that is identifiable at $\bar x$.  Then, for each $r=1,2,3,\ldots$, there exists a point 
$x_r$ in the set $\M \cap B_{1/r}^f(\bar x)$ but outside the set $\M'$.  Consequently we have $x_r \to_f \bar x$ in $\M$, so by definition, the sequence $(x_r)$ is critical, contradicting the identifiability of $\M'$.

Conversely, if the set $\M$ is identifiable but not strongly critical at $\bar x$, then there exists a sequence $x_r \to_f \bar x$ in $\M$ that is not critical.  After taking a subsequence, we can suppose that there exists a value $\delta > 0$ such that $|\nabla f|(x_r) > \delta$ for all $r$.  We claim that the set defined by $\M' = \M \setminus (x_r)$ is identifiable.  If not, there exists a critical sequence $(y_s)$ for $\bar x$ outside $\M'$. Since $\M$ is identifiable, $y_s \in \M$ for all large $s$, and hence $y_s \in (x_r)$, from which we deduce $|\nabla f|(y_s) > \delta$, contradicting the criticality of $(y_s)$.  But now we note that for all $\epsilon > 0$, the point $x_r$ lies in the set $\M \cap B_{\epsilon}^f(\bar x)$ for all large $r$ but lies outside $\M'$, contradicting property (\ref{minimality2}).
\finpf

\noindent
Identifiability and strong criticality play symmetric roles here:  a result analogous to Proposition \ref{minimality} is easy to derive for locally {\em maximal} strongly critical sets.  

We next note an immediate consequence of Proposition \ref{minimality}:  strongly critical identifiable sets (or equivalently, locally minimal identifiable sets), if they exist, must be locally unique.  A subgradient-based version, with stronger assumptions, appeared in \cite[Corollary 4.2]{Hare}.

\begin{prop}[Uniqueness of strongly critical identifiable sets]
Given a function $f$ on a metric space, for any two strongly critical identifiable sets $\M$ and $\M'$ at a point $\bar x$, there exists an $f$-attentive neighborhood $B_{\epsilon}^f(\bar x)$ with 
$\epsilon > 0$ such that
\[
\M \cap B_{\epsilon}^f(\bar x) ~=~ \M' \cap B_{\epsilon}^f(\bar x).
\]
\end{prop}

To illustrate the idea of a strongly critical identifiable set, we consider some simple examples.  At any point $\bar x \in X$, any neighborhood of $\bar x$ is a strongly critical identifiable set for the squared distance function $d^2_{\bar x}$, by Proposition \ref{square}:  the modulus of identifiability is $+\infty$.  On the other hand, if $X$ is a length space, then the set $\M = \{\bar x\}$ is a strongly critical identifiable set for the distance function $d_{\bar x}$, by Theorem \ref{length}:  the modulus of identifiability is $1$.  A simple concrete example to keep in mind is the following.
\begin{exa} \label{basic}
{\rm
On $\R^2$, the function $f(u,v) = |u|+v^2$ has slope
\[
|\nabla f|(u,v) ~=~ 
\left\{
\begin{array}{ll}
\sqrt{1+4v^2} & (u \ne 0) \\
2v & (u=0).
\end{array}
\right.
\]
The set $\M = \{(u,v) : u=0 \}$ is a strongly critical identifiable set at the minimizer zero, and the modulus of identifiability is $1$.  
}
\end{exa}

Unfortunately, as the following example \cite[Example 4.13]{ident} shows, even continuous convex functions do not always admit strongly critical identifiable sets at minimizers.  In such cases, among the collection of sets that eventually capture every critical sequence, no one particular example stands out:  more precisely, there is no locally minimal element. 

\begin{exa}
{\rm
The function 
$f(u,v) = \sqrt{u^2 + v^4}$ has slope zero at its minimizer $(0,0)$, and
\[
|\nabla f|(u,v) ~=~ 
\sqrt{\frac{u^2 + 4v^6}{u^2 + v^4}} \qquad \mbox{for}~ (u,v) \ne (0,0). 
\]
Consider any constant $\alpha > 0$.  Close to $(0,0)$ and outside the set 
\[
\M_{\alpha} ~=~ \{ (u,v) :  |u| \le \alpha v^2 \},
\]
the slope is bounded below by $\frac{\alpha}{\sqrt{1+\alpha^2}} > 0$, so the set $\M_\alpha$ is identifiable at $(0,0)$.  These sets shrink to the set $\M_0$ as $\alpha \downarrow 0$.  However, $\M_0$ is not identifiable, because the critical sequence $(k^{-3},k^{-1})$ for $k=1,2,3,\ldots$ approaches $(0,0)$ from outside $\M_0$.  Thus no locally minimal identifiable set exists at $(0,0)$, and hence, by Proposition~\ref{minimality}, neither does a strongly critical identifiable set.
}
\end{exa} 

The limiting slope gives an equivalent definition of the modulus of identifiability.

\begin{prop}[Identifiability via limiting slope] \label{limiting}
In Definition \ref{def:identifiability}, if the set $\M$ is closed, then the modulus of identifiability is
\bmye \label{limiting-modulus}
\liminf_{x \to_f \, \bar x,~ x \not\in \M} \overline{|\nabla f |}(x).
\emye
\end{prop}

\pf
We just need to show that the modulus (\ref{rate}) is no larger than the quantity (\ref{limiting-modulus}),  the reverse inequality being immediate.
If the claim fails, then there exists a constant $\epsilon > 0$ such that the modulus is at least 
$\overline{|\nabla f |}(x_r) + \epsilon$ for some sequence $x_r \not\in \M$ satisfying $x_r \to_f \bar x$.  Since $\M$ is closed, each distance $\delta_r = \mbox{dist}(x_r,\M)$ is strictly positive, so there exist points $x'_r \in X$ satisfying  
\[
d(x_r,x'_r) < \min\Big\{ \delta_r , \frac{1}{r} \Big\}, \qquad |f(x_r) - f(x'_r)| < \frac{1}{r}, \qquad 
|\nabla f|(x'_r) < \overline{|\nabla f |}(x_r) + \frac{\epsilon}{2}.
\] 
Since $x'_r \not\in \M$ and $x'_r \to_f \bar x$, we have arrived at a contradiction.
\finpf

\noindent
This result may fail if the set $\M$ is not closed.  For example, consider the function $f \colon \R^2 \to \R$ defined by $f(u,v) = \min\{u,0\}$ at the point $\bar x = (0,0)$ in the set $\M$ whose complement is 
$\{(0,v) : v \ne 0\}$.  The modulus (\ref{rate}) is $1$, so $\M$ is identifiable at $\bar x$, but the quantity (\ref{limiting-modulus}) is zero.

As a consequence of Proposition \ref{limiting}, in the special case when the underlying space is Euclidean, the following result shows that our new metric notion of identifiability coincides with the original version \cite[Definition 3.10]{ident}. 

\begin{cor}[Identifiability in Euclidean spaces] \label{original}
Consider a closed function $f$ on $\Rn$ and a point $\bar x$ contained in a closed set $\M \subset \Rn$.  Then 
the set $\M$ is identifiable at $\bar x$ if and only if there exists no sequence $(x_r)$ outside $\M$ satisfying $x_r \to_f \bar x$ and with subgradients $y_r \in \partial f(x_r)$ converging to zero.  
The modulus of identifiability  is
\[
\liminf_{x \to_f \, \bar x,~ x \not\in \M,~ y \in \partial f(x)} |y|.
\]
\end{cor} 

\pf
The formula for the modulus follows from Proposition \ref{limiting} by applying Proposition 
\ref{distance-formula}.  The characterization of identifiability follows. 
\finpf

\noindent
We defer a more detailed investigation of the Euclidean case until Part \ref{part:2}.

\section{Linear growth around identifiable sets} \label{sec:linear}
A crucial geometric feature of objective functions around identifiable sets is often described as ``sharpness''.  An early appearance of this idea was \cite{sharp}, and it was a core ingredient of precursors of identifiability like ${\mathcal {VU}}$-decomposition \cite{MC04} and partial smoothness \cite{Lewis-active}.  Loosely speaking, near a local minimizer, the objective grows in a manner determined entirely by its behavior on any identifiable set, because as we leave this set, the objective grows at a linear rate. 

The original approaches to this linear growth property, including \cite[Proposition 2.10]{Lewis-active}, \cite[Theorem 6.2]{Hare}, and \cite[Theorem D.2]{decay}, rely on smoothness assumptions on the identifiable set $\M$ and the restriction of the objective function $f$ to $\M$.
These assumptions fail in even some simple examples like
\begin{eqnarray*}
f(x) = |x|^{\frac{3}{2}} \quad (x \in \R),  \qquad \bar x = 0, & & \qquad \M = \R \\
f(x) ~=~ x_1^2 + \big| |x_1|^{\frac{3}{2}} - x_2 \big| \qquad (x \in \R^2), 
\qquad \bar x = 0, & &  \qquad \M = \{ x : |x_1|^{\frac{3}{2}} = x_2 \}.
\end{eqnarray*}
By contrast, we develop a simple new linear growth tool for identifiable sets that does not rely on manifolds, subgradients, Euclidean structure, or even nearest-point projections.  Our approach is applicable in any complete metric space, and reveals for the first time the close connection between identifiability and the Ekeland variational principle, which we use in the following form \cite[Theorem 1.4.1]{zal}.

\begin{thm}[Ekeland variational principle  \cite{var_princ}]~
On any complete metric space, if a closed function $f$ is bounded below, then for any value $\epsilon > 0$ and any point $x \in \mbox{\rm dom}\, f$, there exists a point $v$ satisfying 
$f(v) + \epsilon d(v,x) \le f(x)$ and $|\nabla f|(v) \le \epsilon$.
\end{thm}

\begin{thm}[Linear growth] \label{linear}
On a complete metric space, consider a closed function $f$
with slope zero at a point $\bar x$, an identifiable set $\M$ at $\bar x$, and any nonnegative constant $\epsilon$ strictly less than the modulus of identifiability.  Then, for any sequence of points 
$x_r \to_f \bar x$, there exists a sequence of points $v_r  \to_f \bar x$ in $\M$ such that, for all large $r$,
\[
f(v_r) + \epsilon d(v_r,x_r) ~\le~ f(x_r).
\]
\end{thm}

\pf
Denote the modulus of identifiability by $\sigma$.  It suffices to prove the result when 
$0 < \epsilon < \sigma$, since the case $\epsilon = 0$ then follows.

Suppose first that $\bar x$ is a local minimizer.  By redefining $f$ to take the value 
$+\infty$ outside a closed ball centered at the local minimizer $\bar x$, we can assume that 
$\bar x$ is a global minimizer, and so $f$ is bounded below.  Applying the Ekeland principle to each point $x_r$ ensures the existence of points $v_r$ satisfying $|\nabla f|(v_r) \le \epsilon$ and
\[
\epsilon d(v_r,x_r) ~\le~ f(x_r) -  f(v_r) ~\le~ f(x_r) -  f(\bar x) ~\to~ 0.
\]
We deduce $v_r \to \bar x$ and $f(v_r) \to f(\bar x)$, so our assumption about $\epsilon$ ensures 
$v_r \in \M$ for all $r$ larger than some $\bar r$.  Redefining $v_r = \bar x$ for all 
$r \le \bar r$ gives the desired sequence.

Now consider the general case where $|\nabla f|(\bar x) = 0$.  Fix any constant $\delta > 0$ satisfying $2\delta < \sigma-\epsilon$.  The point $\bar x$ locally minimizes the function 
$\tilde f = f + \delta d_{\bar x}$.  Furthermore, since $f = \tilde f + (-\delta d_{\bar x})$, and the function $-\delta d_{\bar x}$ is $\delta$-Lipschitz, Proposition~\ref{slope-tools} implies, for all points $x$, the inequality $|\nabla\tilde f|(x) \ge |\nabla f|(x) - \delta$.  Taking the $\liminf$ of both sides as $x \to_f \bar x$ (or equivalently as $x \to_{\tilde f} \bar x$) outside the set $\M$ shows that $\M$ is also identifiable for the function~$\tilde f$, with modulus at least
$\sigma - \delta > \epsilon + \delta$.

We can now apply the result proved by the first argument, using the function $\tilde f$ in place of the function $f$ and constant $\epsilon + \delta$ in place of $\epsilon$.  We deduce that there exists a sequence of points $v_r  \to \bar x$ in $\M$ satisfying 
$\tilde f(v_r) \to \tilde f(\bar x)$ and 
\[
\tilde f(v_r) + (\epsilon + \delta) d(v_r,x_r) ~\le~ \tilde f(x_r) \qquad \mbox{for all large}~ r.
\]
Consequently $f(v_r) \to f(\bar x)$, and for all large $r$ we have
\[
f(v_r) + \delta d(v_r,\bar x) + (\epsilon + \delta) d(v_r,x_r) ~\le~ f(x_r) + \delta d(x_r,\bar x).
\]
Our conclusion follows by the triangle inequality.
\finpf

\noindent
The slope zero assumption in Theorem \ref{linear} cannot be relaxed to criticality. For example, the function $f(x) = \min\{0, x\}$ for $x \in \R$ has a critical point $\bar{x} = 0$ for which the set $\M = \R_+$ is identifiable, but the result fails for the sequence $x_r = -\frac{1}{r}$.

We next derive simple optimality conditions from Theorem \ref{linear} (Linear growth).  We model these  on \cite[Theorem 6.2]{Hare} and \cite[Proposition 7.2 and Corollary 7.3]{ident}, but replacing the  subdifferential arguments for those Euclidean results with our simpler Ekeland-based approach, valid in an arbitrary complete metric space.  In the case of max functions (Section \ref{max}), the results reduce to classical active-set properties.

\begin{cor}[Sufficient condition for optimality]~
On a complete metric space, suppose that a closed function $f$ has slope zero at a point $\bar x$, and consider any identifiable set $\M$ at $\bar x$.  Then $\bar x$ is a local minimizer if and only if it is a local minimizer relative to $\M$.  Furthermore, $\bar x$ is a strict local minimizer if and only if it is a strict local minimizer relative to $\M$.
\end{cor}

\pf
To prove the first equivalence,
suppose that $\bar x$ is not a local minimizer, so there exists a sequence $x_r \to \bar x$ such that $f(x_r) < f(\bar x)$.  Closedness implies $f(x_r) \to f(\bar x)$, so using Theorem \ref{linear} (Linear growth), we deduce the existence of a constant $\epsilon > 0$ and a sequence of points $v_r  \to \bar x$ in $\M$ satisfying, for all large $r$,
\[
f(v_r) ~\le~ f(v_r) + \epsilon d(v_r,x_r) ~\le~ f(x_r) ~<~ f(\bar x),
\]
so $\bar x$ is not a local minimizer relative to $\M$.  The converse is trivial.  The second equivalence concerns strict local minimizers of $f$, namely points $\bar x$ such that $f(x) > f(\bar x)$ for all points $x \ne \bar x$ near $\bar x$.  The proof is very similar.
\finpf

As a further illustration along classical lines, we show that quadratic growth rates are determined by the growth rate on any identifiable set.

\begin{cor}[Quadratic growth]
On a complete metric space, suppose that a closed function $f$ has slope zero at a point $\bar x$, and consider any identifiable set 
$\M$ at $\bar x$.  Then $f$ has quadratic growth around $\bar x$ if and only if it has quadratic growth around 
$\bar x$ relative to $\M$.  Indeed, the two growth rates are identical:
\[
\liminf_{x \to \bar x \atop x \ne \bar x} \frac{f(x) - f(\bar x)}{d(x,\bar x)^2}
~=~
\liminf_{x \to \bar x,~ x \in \M \atop x \ne \bar x} \frac{f(x) - f(\bar x)}{d(x,\bar x)^2}.
\]
\end{cor}

\pf
Denote the left-hand and right-hand sides by $\alpha$ and $\beta$ respectively.  Clearly 
$\alpha \le \beta$.  Suppose in fact $\alpha < \beta$.  Choose any value $\gamma$ in the interval $(\alpha,\beta)$.  Since $\gamma > \alpha$
there is a sequence of points $x_r \to \bar x$ satisfying $f(x_r) - f(\bar x) < \gamma d(x_r,\bar x)^2$.
Taking the $\limsup$ of both sides shows $\limsup_r f(x_r) \le f(\bar x)$.  Since $f$ is closed, we also know $\liminf_r f(x_r) \ge f(\bar x)$, so we deduce $f(x_r) \to f(\bar x)$.  By Theorem \ref{linear} (Linear growth), there exists a constant $\epsilon > 0$ and a sequence of points $v_r \to \bar x$ in $\M$ satisfying $f(v_r) + \epsilon d(v_r,x_r) \le f(x_r)$ for all large $r$.  Hence $v_r \ne \bar x$ for all large $r$, since otherwise we arrive at the contradiction $f(\bar x) + \epsilon d(\bar x,x_r) < f(\bar x) + \gamma d(x_r,\bar x)^2$.
Since $\gamma < \beta$, for all large $r$ we must have 
\begin{eqnarray*}
\gamma d(v_r,\bar x)^2  
&<&
 f(v_r) - f(\bar x) ~\le~ f(x_r) - \epsilon d(v_r,x_r) - f(\bar x) \\
&<& 
\gamma d(x_r,\bar x)^2  - \epsilon d(v_r,x_r) ~\le~ \gamma \big( d(v_r,\bar x) + d(v_r,x_r) \big)^2  - \epsilon d(v_r,x_r).
\end{eqnarray*}
We deduce the inequality $0 < d(v_r,x_r) \big(2d(v_r,\bar x) + d(v_r,x_r) - \epsilon \big)$,
which is impossible for $r$ sufficiently large.
\finpf

\noindent
In Section \ref{KL}, we prove an analogous result for a growth condition fundamental to convergence analysis for optimization algorithms.

Theorem \ref{linear} asserts only the existence of a ``shadow'' sequence $v_r$. If the underlying space is Euclidean and $\M$ is a smooth manifold restricted to which $f$ is smooth, then \cite[Theorem D.2]{decay} is more geometrically descriptive:  for small $\epsilon$ we can take as $v_r$ the nearest-point projection $\mathrm{proj}_{\mathcal{M}} (x_r)$.  However, we cannot expect this more generally.  First, unless the metric space $X$ is proper, the requisite nearest points may not exist.  Secondly, the argument of \cite[Theorem D.2]{decay} relies on smoothness of $\M$ and the corresponding restriction of $f$, and furthermore on a characterization of identifiable manifolds using subgradients \cite{ident2} that is very Euclidean in spirit.  

However, by relaxing our notion of projection, we can leverage Theorem \ref{linear} (Linear growth) to prove an analogous result in great generality.  This new result uses {\em approximate} projections of points $x$ onto $\M$:  after fixing a constant $C \ge 1$, we consider points $v \in \M$ such that 
$d(v, x) \le C \mbox{dist}(x, \M)$.

We need some standard metric space notions.  As usual, we consider a metric space $(X,d)$, a subset $\M$, and a function $f \colon X \to (-\infty,+\infty]$.  By restricting the original metric $d$ to 
$\M$, we arrive at the {\em induced} metric on $\M$, also denoted $d$ for simplicity.  On the resulting metric space $(\M,d)$, we denote the {\em restriction} of $f$ by $f|_\M$. In $\M$, if any two points can be joined by a rectifiable curve in $\M$, then we can also define the {\em intrinsic} metric $d^\M$:  the intrinsic distance between any two points is the infimum of the lengths of curves joining them.  If there exists a constant $\mu>0$ such that $d^\M(x,y) \le \mu d(x,y)$ for all points $x,y \in \M$, then we say that $\M$ is {\em normally embedded} \cite{birbrair-mostowski}.  (By the triangle inequality, 
$d(x,y) \le d^\M(x,y)$ always holds.) Clearly the property is preserved by bi-Lipschitz homeomorphisms.  Examples include convex sets and compact submanifolds in $\Rn$.  Curves in normally embedded sets have the same length with respect to both the induced and intrinsic metrics \cite[Proposition 2.3.12]{burago}, and furthermore, when rectifiable, can be reparametrized to have unit speed (with respect to both metrics) \cite[Proposition~2.5.9]{burago}.   

\begin{thm}[Linear growth from approximate projection] \label{linear-proj}
On a complete metric space, consider a closed function $f$ with slope zero at a point $\bar x$, and a normally embedded identifiable set $\M$ at $\bar x$ such that the restriction $f|_\M$ is continuous.
\begin{itemize}
\item
If $\M$ is strongly critical at $\bar x$, then the slope $| \nabla(f|_\M )|$ is continuous at 
$\bar x$.
\item
Conversely, if $|\nabla (f|_\M)|$ is continuous at $\bar x$, then
for any constant $C \ge 1$, there exists a constant 
$\epsilon > 0$ such that all points $x$ close to $\bar x$ with value $f(x)$ close to $f(\bar x)$ satisfy the property
\[
v \in \M \quad \mbox{and} \quad d(v, x) \le C \mbox{\rm dist}(x, \M) \quad \Rightarrow \quad f(v) + \epsilon d(v, x) \le f(x).
\]
If, furthermore, $\M$ is locally closed, then $\M$ is strongly critical at $\bar x$.
\end{itemize}
\end{thm}

\pf
The first claim is clear, because $|\nabla f|(x) \ge | \nabla(f|_\M )|(x)$ for all points $x \in \M$.
Turning to the converse direction, by way of contradiction, consider sequences of points $x_r \to_f \bar x$ and $v_r$ in $\M$ satisfying the inequalities $d(v_r, x_r) \le C \mbox{dist}(x_r, \M)$ and $f(v_r) + \frac{1}{r} d(v_r,x_r) > f(x_r)$ for all $r$.
By Theorem \ref{linear} (Linear growth), there exists $\epsilon > 0$ and a sequence $w_r \rightarrow \bar x$ in $\M$ satisfying $f(w_r) + \epsilon d(w_r,x_r) \le f(x_r)$ for all large $r$.
From the inequalities $d(v_r, x_r) \le C \mbox{dist}(x_r, \M) \le C d(w_r, x_r)$,
we deduce
\[
d(w_r, v_r) ~\le~ d(w_r, x_r) + d(x_r, v_r) ~\le~ (1 + C) d(w_r, x_r).
\]
Combining the inequalities shows, for all large $r$,
\begin{eqnarray*}
f(v_r) - f(w_r) &>&
\epsilon d(w_r,x_r) - \frac{1}{r} d(v_r,x_r) ~\ge~ \big(\epsilon - \frac{C}{r}\big) d(w_r,x_r) 
~\ge~ 
\frac{\epsilon}{2} d(w_r,x_r) \\
& \ge & \frac{\epsilon}{2(1 + C)} d(w_r,v_r).
\end{eqnarray*}

Since the set $\M$ is normally embedded, there exists a number $\mu > 0$ such that the intrinsic metric $d^\M$ satisfies $d^\M(w,v) \le \mu d(w,v)$ for all points $w,v \in \M$.
In particular, since $w_r,v_r \in \M$, we have
\[
d^\M(w_r,v_r) ~\le~ \mu d(w_r,v_r) ~\le~ 
\frac{1}{\delta} \big(f(v_r) - f(w_r)\big)
\]
where $\delta = \frac{\epsilon}{2\mu(1+C)}$, so there exist unit-speed curves 
$\gamma_r: [0, l_r] \to \M$ from $w_r$ to $v_r$ with length 
$l_r \le 2 d^\M(w_r,v_r) \le \frac{2}{\delta}\big(f(v_r) - f(w_r)\big)$.  The corresponding continuous functions
$g_r = f \circ \gamma_r$ satisfy $g_r (l_r) - g_r(0) > \frac{\delta}{2} l_r$.
Let $t_r$ be the smallest number $t \in [0, l_r]$ such that
\[
g_r(s) ~\ge~ g_r(l_r) - \frac{\delta}{2} (l_r-s)  \qquad  \mbox{for all}~ s \in [t, l_r].
\]
Continuity implies $t_r > 0$ and $g_r(t_r) = g_r(l_r) - \frac{\delta}{2} (l_r-t_r)$. Moreover, for each large $r$, there exists a sequence of points $s^r_k \nearrow t_r$ such that 
$g_r(s^r_k) < g_r(l_r) - \frac{\delta}{2}(l_r - s^r_k)$ for all $k = 1,2,3,\ldots$.
We deduce 
\[
\limsup_{t \nearrow t_r} \frac{g_r(t_r) - g_r(t)}{t_r - t} ~\ge~ \limsup_{k \to \infty} \frac{g_r(t_r) - g_r(s^r_k)}{t_r - s^r_k} ~\ge~ \frac{\delta}{2}.
\]
Consider the points $y_r = \gamma_r (t_r) \in \M$. Since $\gamma_r$ is unit-speed, we have
\begin{eqnarray*}
|\nabla(f|_\M)|(y_r)
&=&
\limsup_{x \tiny{\xrightarrow{\M}} y_r,\, x \neq y_r} \frac{f(y_r) - f(x)}{d(y_r, x)} 
~\ge~  
\limsup_{x \tiny{\xrightarrow{\M}} y_r,\, x \neq y_r} \frac{f(y_r) - f(x)}{d^{\M} (y_r, x)} \\
&\ge& 
\limsup_{t \nearrow t_r} \frac{f(y_r) - f\big(\gamma_r (t)\big)}{d^{\M} \big(\gamma_r(t_r), \gamma_r (t)\big)} ~\ge~ 
\limsup_{t \nearrow t_r} \frac{g_r(t_r) - g_r(t)}{t_r-t} ~\ge~ \frac{\delta}{2}.
\end{eqnarray*}
Furthermore we have
\[
d(y_r, \bar x) ~\le~ d^{\M} (y_r, w_r)+ d(w_r, \bar x) ~\le~ 2 d^{\M} (v_r, w_r) + d(w_r, \bar x) \to 0,
\]
contradicting the continuity of the slope $|\nabla(f|_\M)|$ at   $\bar x$.

It remains to prove that the set $\M$ is strongly critical at the point $\bar x$.  Fix any constant 
$\delta > 0$.  By the slope continuity assumption, there exists a constant $\epsilon > 0$ such that 
$|\nabla(f|_\M)|(x) < \delta$ for all points $x$ in the ball $B_{2\epsilon}(\bar x)$.    Since $\M$ is locally closed, we can furthermore assume that such $x$ have approximate projections $v \in \M$ with $d(v, x) \le 2\mathrm{dist}(x, \M)$.
Our argument so far shows that, after shrinking $\epsilon$, we can assume that these approximate projections satisfy $f(v) \le f(x)$.  

We claim $|\nabla f|(x) \le 3\delta$ for all points $x \in \M \cap B_{\epsilon}(\bar x)$. We can assume that $x$ is not a local minimizer for $f$, since otherwise the claim is immediate.  Choose a sequence $x_r \to x$ in $B_{\epsilon}(\bar x)$ such that 
\[
q_r ~=~ \frac{f(x)-f(x_r)}{d(x,x_r)} ~\to~ |\nabla f|(x)
\]
and choose approximate projections $v_r \in \M$ with $d(v_r, x_r) \le 2\mathrm{dist}(x_r, \M)$.  By definition, these points satisfy
$d(v_r,x_r) \le 2 d(x,x_r)$ and hence 
\[
d(x,v_r) ~=~ \big(d(x,v_r) - d(v_r,x_r)\big) + d(v_r,x_r) ~\le~ 3d(x,x_r),
\]
implying in particular $v_r \to x$.  Moreover, as we have seen, $f(v_r) \le f(x_r)$ for all large $r$.  Since $|\nabla(f|_\M)|(x) < \delta$, we know $f(x) - f(v_r) < \delta d(x,v_r)$ for all large $r$, and hence
\[
q_r ~<~ \frac{\delta d(x,v_r)}{d(x,x_r)} 
~\le~ 3\delta.
\]
Letting $r \to \infty$ proves the claim.  Since $\delta$ was arbitrary, the result follows.
\finpf

\noindent
While we have assumed, for simplicity that the set $\M$ is normally embedded, we only need that assumption to hold {\em locally} around the point $\bar x$, in the obvious sense.

In the most basic case, the underlying space is Euclidean and the identifiable set $\M$ is a smooth submanifold, so the nearest point projection mapping
$\mathrm{proj}_{\mathcal{M}}$ is single-valued near any point $\bar x \in \M$.  Setting 
$v=\mathrm{proj}_{\mathcal{M}}(x)$ in Theorem~\ref{linear-proj}, we then recover \cite[Theorem D.2]{decay} as a special case.  However, this new result applies more broadly, including situations when $\M$ is not a manifold, as we illustrate later.

\section{The Kurdyka-{\L}ojasiewicz property} \label{KL}
Our initial motivation for the preceding section was the idea of sharpness, a notion highlighted in the optimization literature in \cite{sharp}.  We now return to that idea in order to discuss the nonsmooth Kurdyka-{\L}ojasiewicz inequality \cite{loja,Lewis-Clarke}, a property that has had far-reaching impact on complexity analysis for optimization.  We begin with the simple idea of sharpness at a point.

\begin{defn}
{\em
On a metric space $(X,d)$, a point $\bar x$ is a {\em sharp local minimizer} of a function $f \colon X \to (-\infty,+\infty]$ if there exists a constant $\epsilon > 0$ such that $\bar x$ locally minimizes 
$f - \epsilon d_{\bar x}$.
}
\end{defn}

\noindent
For example, any point $\bar x \in X$ is a sharp local minimizer of the distance function $d_{\bar x}$.  Notice that any sharp local minimizer is also a strict local minimizer.  An immediate consequence of  Theorem \ref{linear} (Linear growth) is the following result.

\begin{cor}[Sharp minimizers] \label{cor:sharp}
On a complete metric space, if a closed function $f$ has slope zero at a point $\bar x$, and the set 
$\{\bar x\}$ is identifiable,  then $\bar x$ is a sharp local minimizer.
\end{cor}

\noindent
The following example shows that the converse may fail.  

\begin{exa}[Sharp minimizers may not be identifiable] \label{unidentifiable}
{\rm
The locally Lipschitz function defined on $\R$ by $f(0) = 0$ and $f(x) = |x| + x^2 \sin \frac{1}{x}$ for $x \ne 0$ satisfies $f(x) \ge \frac{1}{2}|x|$ for all small $x$, and hence has a sharp minimizer at $0$. However, $f$ has derivative zero at each point $\frac{1}{2\pi r}$ for $r=1,2,3,\ldots$.  These points approach zero, so the set $\{0\}$ is not identifiable.
}
\end{exa}

To handle the more practical situation where minimizers may not be strict, \cite{sharp} introduced the following weaker property, which we use in a local form.

\begin{defn}
{\rm
On a metric space $X$, a point $\bar x$ is a {\em weak sharp minimizer} (locally) for a function 
$f \colon X \to (-\infty,+\infty]$ if there exists a constant $\epsilon > 0$ such that, for
the level set ${\mathcal L} = \{x : f(x) \le f(\bar x)\}$,
the point $\bar x$ locally minimizes the function $f - \epsilon \mbox{dist}(\cdot,{\mathcal L})$.
}
\end{defn}

\noindent
By analogy with Corollary \ref{cor:sharp}, we now draw a new relationship between weak sharpness and identifiability.

\begin{cor}[Weak sharp minimizers] \label{weak}
On a complete metric space, consider a closed function $f$ and a local minimizer $\bar x$.  If the level set of points $x$ with value  $f(x) \le f(\bar x)$ is identifiable at $\bar x$ for $f$, then $\bar x$ is a weak sharp minimizer.
\end{cor}

\pf
Denote the level set by ${\mathcal L}$.
Suppose $\epsilon > 0$ is strictly less than the modulus of identifiability of ${\mathcal L}$ for the function $f$ at the point $\bar x$:
\[
\liminf_{x \to_f \, \bar x, ~f(x) > f(\bar x)} |\nabla f |(x) ~>~ \epsilon.
\]
If the result fails, then there exists a sequence of points $x_r \to \bar x$ in $X$ satisfying 
\bmye \label{strict}
f(x_r) ~<~ f(\bar x) + \epsilon \cdot \mbox{dist}(x_r,{\mathcal L}).
\emye
Since $\bar x$ is a local minimizer, for large $r$ we have $f(x_r) \ge f(\bar x)$ and so inequality (\ref{strict}) implies 
$x_r \not\in {\mathcal L}$. The right-hand side converges to $f(\bar x)$, whence so does $f(x_r)$.
By Theorem \ref{linear} (Linear growth), there exists a sequence of points $v_r  \to \bar x$ in ${\mathcal L}$ satisfying, for large $r$, the inequality $f(v_r) + \epsilon d(v_r,x_r) \le f(x_r)$.
For large $r$ we know $f(v_r) \ge f(\bar x)$, so $\epsilon d(v_r,x_r) \le  f(x_r) - f(\bar x)$,
contradicting inequality (\ref{strict}).
\finpf

\noindent
Example \ref{unidentifiable} shows that the converse of Corollary \ref{weak} may fail.

The Kurdyka-{\L}ojasiewicz property may be thought of as a modification of the assumption of Corollary~\ref{weak}, obtained by first shifting and truncating the function $f$, replacing each value $f(x)$ by $[f(x) - f(\bar x)]^+$, where $u^+$ denotes the positive part of the value $u$, and then rescaling those new function values using the following standard tool.

\begin{defn} 
{\rm
A {\em desingularizer} is a function $\phi \colon [0,+\infty] \to [0,+\infty]$ that is continuous on $[0,+\infty)$, satisfies $\phi(0) = 0$ and $\phi(\tau) \to \phi(+\infty)$ as $\tau \to +\infty$, and has continuous strictly positive derivative on $(0,+\infty)$.
}
\end{defn}

\noindent
Typically we use desingularizers of the form $\phi(\tau) = \tau^{1-\alpha}$, where the {\em KL~exponent} $\alpha$ lies in the interval $[0,1)$.  Such desingularizers, for example, suffice for the original {\L}ojasiewicz proof of the KL property for real-analytic functions on Euclidean spaces \cite{loj}, and hence also for real-analytic functions on analytic manifolds \cite[Section 9]{kurdyka-mostowski-parusinski}.  Unlike some definitions, we do not explicitly require concavity.

We now present a definition of the Kurdyka-{\L}ojasiewicz property in metric spaces.  Our definition is inspired by a version of the KL property from \cite[Corollary 4]{tailwag}, reframed in the language of identifiability.  In the Euclidean case, some earlier subgradient-based versions of the KL property \cite{attouch-projections} reduce precisely to this metric definition;  the original nonsmooth property \cite{loja,Lewis-Clarke} (again subgradient-based) is slightly stronger, as we shortly discuss.  Researchers  on Riemannian nonsmooth optimization have proposed generalizations to manifolds \cite{bento,Hosseini:2015}, but by comparison with the subgradient approach in general, the metric definition we present is simple and direct.  By contrast with \cite[Section 2]{tailwag}, which aimed to illuminate the KL inequality as a metric regularity property, the new connection to identifiability that we present allows us to deploy the linear growth results of Section~\ref{sec:linear}.

\begin{defn} \label{def:kl}
{\rm
On a metric space $X$, a function $f \colon X \to (-\infty,+\infty]$ satisfies the {\em KL property} at a point $\bar x \in \mbox{dom}\, f$ if there exists a desingularizer $\phi$ such that the level set of points $x$ satisfying $f(x) \le f(\bar x)$ is identifiable at $\bar x$ for the composite function $g$ defined on $X$ by $g(x) = \phi\big([f(x)-f(\bar x)]^+\big)$.
}
\end{defn}

We make some immediate comments.  The level set for $f$ in the definition is also the level set defined by $g(x) \le 0$, for the composite function $g$.  Hence, by Corollary \ref{weak}, when the metric space $X$ is complete and $f$ is closed, the KL property implies that the point $\bar x$ is a weak sharp minimizer of the composite function $g$.  

For closed functions $f$, the KL property amounts to the existence of $\delta > 0$ such that
$\big|\nabla\big(\phi(f(\cdot)-f(\bar x))\big)\big|(x) \ge \delta$
for all points $x$ near $\bar x$ with value $f(x)$ near and strictly larger than $f(\bar x)$, as stated in \cite[Corollary 4]{tailwag}.  Using a simple chain rule \cite[Lemma 4.1]{aze2017}, we can rewrite the left-hand side as
$\phi'\big(f(x)-f(\bar x)\big) \cdot |\nabla f|(x)$.  
When the underlying space $X$ is Euclidean, Proposition \ref{distance-formula} transforms our definition into an equivalent but more familiar inequality:
\bmye \label{weak-KL}
 \phi'\big(f(x)-f(\bar x)\big) \cdot \mbox{dist}\big(0,\partial f(x)\big) ~\ge~ \delta.
\emye
 The original Euclidean nonsmooth property \cite{Lewis-Clarke} (sometimes called the {\em strong} KL inequality \cite{noll-kl}), while theoretically stronger, since it replaces the subdifferential $\partial f(x)$ in inequality (\ref{weak-KL}) by the larger {\em Clarke} subdifferential, in practice seems no stronger:  for more discussion, see \cite{lewis-tian}.  In Euclidean spaces, the strong KL property always holds for functions $f$ that are semi-algebraic, or, more generally, tame \cite{Lewis-Clarke}.
 
Remarkably, by-now-standard KL-based convergence analyses like \cite[Theorem 24]{tailwag} often transpire on inspection to be entirely metric in nature.  This new perspective clarifies proofs and extends their reach to settings such as manifolds and more general metric spaces.  Before illustrating this view, we require a standard definition.

\begin{defn} 
{\rm
A function $f$ on a metric space is {\em continuous on slope-bounded sets} \cite{ambrosio} if sequences of points $x_r \to x$ with values $f(x_r)$ and slopes $|\nabla f|(x_r)$ uniformly bounded must satisfy $f(x_r) \to f(x)$. 
}
\end{defn}

\noindent
In particular, on Euclidean spaces, all closed convex functions are continuous on slope-bounded sets, as, more generally, are all closed functions that are subdifferentially continuous in the sense of \cite{VA}.

We next illustrate KL-based analysis in metric spaces by proving convergence of the proximal point method.  The argument follows that for \cite[Theorem 24]{tailwag}. In that result, the setting is a Hilbert space and the objective must be weakly convex.  However, on inspection, as we show, the result has an entirely metric version.  For comparison, \cite{bacak} also proves convergence of the proximal point method in metric spaces, but requires the space to be of nonpositive curvature and the objective to be geodesically convex;  furthermore, the argument only guarantees a notion of convergence weaker than the convergence in distance that we prove next.
 
\begin{exa}[Proximal sequences in metric spaces]
{\rm
On a complete metric space, we consider a closed function 
$f$ that is bounded below, and continuous on slope-bounded sets.  We assume that there exists a value 
$\rho > \inf f$ and a desingularizer $\phi$ such that
\[
\inf f < f(x) < \rho \qquad \Rightarrow \qquad \big|\nabla \big(\phi(f(\cdot) - \inf f)\big)\big|(x) ~\ge~ 1.
\]
A simple open cover argument (cf.\ \cite[Lemma 1(ii)]{attouch-bolte}) shows that this uniform version of the KL property holds if the level set $\{ x : f(x) \le \rho\}$ is compact and the KL property holds at every minimizer.  We consider a sequence of points $(x_k)$ arising from the proximal point iteration:
\[
x_{k+1} ~\in~ \mbox{prox}_f(x_k) \qquad (k=0,1,2,\ldots).
\]
As we discussed in Section \ref{sec:Identifiability}, in quite general settings there exist such sequences starting from any initial point $x_0 \in \mbox{dom}\, f$.  The values $f(x_k)$ are nonincreasing and bounded below, and Lemma \ref{lem:proximal} implies $|\nabla f|(x_{k+1}) \le d(x_{k+1},x_k)$ for all $k$.  We furthermore assume the initial value satisfies $f(x_0) < \rho$.
We claim that the iterates $x_k$ must converge to a critical point.  We lose no generality in supposing $f(x_k) > \inf f$ for all $k$.  Consider the level sets $L_k = \{x : f(x) \leq f(x_k)\}$.
By definition, all points $x \in X$ satisfy \
\[
f(x_{k+1}) + \frac{1}{2}d^2(x_k,x_{k+1}) ~\le~ f(x) + \frac{1}{2}d^2(x_k,x),
\]
so all $x$ in the level set $L_{k+1}$ satisfy $d(x_k,x) \ge d(x_k,x_{k+1})$.
We deduce the equation $d(x_k, x_{k+1}) = \mathrm{dist}(x_k, L_{k+1})$.  The right-hand side is no larger than the Hausdorff distance between the sets $L_k$ and $L_{k+1}$, which by \cite[Corollary~4]{tailwag} is bounded above by $\phi\big(f(x_k)\big) - \phi\big(f(x_{k+1})\big)$.  Thus trajectory $\sum_k d(x_k,x_{k+1})$ has finite length, so by completeness, $x_k$ converges to some point $x^*$ and furthermore $|\nabla f|(x_k) \to 0$.  Continuity on slope bounded sets implies $f(x_k) \to f(x^*)$, so $x^*$ is a critical point.
}
\end{exa}

We are now ready for one of our main new results.  In the previous section, we used Theorem \ref{linear} (Linear growth) to show how a function, around any point where its slope is zero, inherits various growth properties from the corresponding properties restricted to any identifiable set.  We now extend that list to include the KL property.  The possibility of such an extension is perhaps not surprising:  in general metric spaces, the KL inequality is closely related to a growth property \cite[Corollary 4]{tailwag}, and in the convex Euclidean case the ideas are essentially equivalent \cite[Section 3]{bolte_complex_kur}.  However, this powerful technique for recognizing the KL property through identifiability is new even in the Euclidean case.  Furthermore, we show that the unrestricted property holds with the same desingularizer (and, in particular, KL exponent) as the restricted property.

Recall that on any subset $\M$ of the metric space $(X,d)$, we can consider the induced metric, where the distance between points $x,y \in \M$ is just $d(x,y)$. The set $\M$ equipped with the induced metric then forms a metric space $(\M, d)$.

\begin{thm}[Identifiability and the KL property]~ \label{kl-identifiable}
On a complete metric space $(X,d)$, suppose that a closed function $f$ has slope zero at a point $\bar x$.  Consider any identifiable set $\M$ at $\bar x$, and a desingularizer $\phi$ that is concave (or, more generally, that satisfies $\liminf_{\tau \searrow 0} \phi'(\tau) > 0$). Then the following properties are equivalent.
\begin{enumerate}
\item[{\rm (a)}]
The function $f$ has the KL property at $\bar x$ with desingularizer $\phi$.
\item[{\rm (b)}]
On the metric space $(\M, d)$, the function $f|_\M$ has the KL property at $\bar x$ with desingularizer $\phi$.
\item[{\rm (c)}]
The function $f + \delta_\M$ has the KL property at $\bar x$ with desingularizer $\phi$.
\end{enumerate}
\end{thm}

\pf
Properties (b) and (c) differ only in notation.  Observe that identifiability ensures the existence of a constant $\epsilon > 0$ such that 
$|\nabla f|(x) > 2\epsilon$ for all points $x \not\in \M$ near $\bar x$ with value $f(x)$ near 
$f(\bar x)$. Without loss of generality assume $f(\bar x) = 0$.

Suppose that property (b) holds.  Then there exists a constant $\delta> 0$ such that
$\big|\nabla \big(\phi \circ f|_\M\big)\big|(x) \ge \delta$
for all points $x \in \M$ near $\bar x$ satisfying $0 < f(x) < \delta$.
From this inequality, the definition of slope then implies 
\bmye \label{inside}
\big|\nabla \big(\phi \circ f \big)\big|(x) ~\ge~ \delta,
\emye
since the left-hand side is no smaller than the previous left-hand side.

On the other hand, consider points $x \not\in \M$ near $\bar x$ satisfying $0 < f(x) < \delta$. 
Shrinking $\delta$ if necessary, we can ensure $\phi'\big(f(x)\big) \ge \delta$ and 
 $|\nabla f|(x) > 2\epsilon$. The chain rule implies 
 $\big|\nabla \big(\phi \circ f \big)\big|(x) \ge \delta \cdot 2\epsilon > 0$.  Property (a) follows from 
inequality (\ref{inside}).

Conversely, if property (a) holds, then there exists a constant $\delta > 0$ such that both the inequalities $\big|\nabla \big(\phi \circ f \big)\big|(x) > \delta$ and $\phi'\big(t\big) \ge \delta$ hold for all points $x$ satisfying $d(x,\bar x) < \delta$ and $0 < f(x) < \delta$, and all
$t \in (0, \delta)$.  Shrinking $\delta$ if necessary, identifiability ensures that all such points $x$ outside $\M$ also satisfy $|\nabla f|(x) > 2\epsilon$.  On the other hand, for any such point $x$ in $\M$, there exists a sequence of points $x_r \to x$, each satisfying 
$\phi\big(f(x)\big) - \phi\big(f(x_r)\big) > \delta \cdot d(x,x_r)$.
By Theorem \ref{linear} (Linear growth), there exist points $v_r \in \M$, satisfying
$f(v_r) + \epsilon d(v_r,x_r) \le f(x_r)$
for all large $r$, so
\begin{eqnarray*}
\phi\big(f(x)\big) - \phi\big(f(v_r)\big) 
& > & 
\delta \cdot d(x,x_r) + \phi\big(f(x_r)\big) - \phi\big(f(v_r)\big) \\
& = & 
\delta \cdot d(x,x_r) + \phi'(t_r)\big( f(x_r) - f(v_r) \big)
\end{eqnarray*}
for some value $t_r \in [f(v_r),f(x_r)]$, by the mean value theorem. Hence we deduce
\[
\phi\big(f(x)\big) - \phi\big(f(v_r)\big) ~>~  \delta \cdot d(x,x_r) ~+~ \delta \cdot  \epsilon  \cdot d(v_r,x_r) ~\ge~ \gamma \cdot d(x,v_r),
\]
where $\gamma = \delta \cdot \min \{1 , \epsilon \} > 0$, so 
$|\nabla(\phi \circ f|_\M)|(x) \ge \gamma$.
Property (b) follows.
\finpf

\begin{rem}
{\rm
By restricting attention to desingularizers $\phi$ satisfying the condition \mbox{$\liminf_{\tau \searrow 0} \phi'(\tau) > 0$}, we lose no essential generality. If a desingularizer $\phi$ fails this condition, then the function defined by $\varphi (t) = \int_0^t \max\{1, \phi'(\tau)\} d\tau, \forall\, t \geq 0$, is another desingularizer, and it maintains the KL property and satisfies the condition.
}
\end{rem}

\subsection*{Talweg curves and identifiability}
In \cite{tailwag}, the authors explore characterizations of the KL property via the finite length of ``talweg'' curves (German for ``valley paths''). For a suitably isolated critical point $\bar{x}$, given any constant 
$R > 1$, they consider the mapping 
$\chi \colon \big(f(\bar x),+\infty\big) \tto X$
defined, for a closed bounded neighborhood $D$ of $\bar{x}$, by
\[
\chi(r) ~=~ \Big\{
x \in D : f(x) = r, ~\overline{|\nabla f|} (x) \leq 
R \inf_{y \in D \atop {f(y)=r}} \overline{|\nabla f|} (y)
\Big\}.
\]
A {\em talweg} is a curve 
$x \colon \big(f(\bar x), \rho\big) \to X$ satisfying $x(r) \in \chi(r)$ for all $r $.

Consider, for example, Example \ref{basic}, where $f(u,v) = |u| + v^2$.  The unique critical point 
$\bar x = (0,0)$ is also the strict global minimizer, with value $f(\bar x) = 0$.  Let $D$ be the unit disk. Suppose $0<r<\frac{1}{4R^2}$.  Then we have
\[
R \inf_{y \in D \atop {f(y)=r}} \overline{|\nabla f|} (y) ~=~ 2R\sqrt{r} ~<~1,
\]
so $\chi(r) = \{(0,\pm\sqrt{r})\}$.  Thus any talweg curve $x: (0, r) \to \R^2$ lies in the identifiable set $\M = \{0\} \times \R$ for all small $r > 0$.

We can generalize this observation using the following assumption.

\begin{ass}
{\rm
At some strict local minimizer $\bar x$, the set $\M$ is strongly critical, identifiable, and locally connected. Furthermore, $\bar x$ is not an isolated point of $\M$, and the restriction $f|_\M$ is continuous. 
}
\end{ass}

\noindent
With this assumption, we claim
\bmye \label{small-slope}
\inf_{y \in D \atop {f(y)=r}}  \overline{|\nabla f|} (y) ~\to~ 0 \quad \mbox{as}~ r \downarrow f(\bar x).
\emye
If not, then there exists a value $\delta > 0$ and a sequence $f(\bar x) < r_i \downarrow f(\bar x)$ such that 
\bmye \label{big-slope}
\overline{|\nabla f|} (y) ~>~ \delta \quad \mbox{for all}~ y \in D ~\mbox{such that}~ f(y) = r_i~ \mbox{for some}~ i.
\emye
We can assume the existence of points $y_i \in \M$ approaching $\bar x$ and satisfying $f(y_i) = r_i$ for all $i$.  Otherwise, taking a subsequence, there would exist an open neighborhood $U$ of $\bar x$ such that $f(x) \ne r_i$ for all $x \in U$ and all $i$.  We can assume, after shrinking, that $U$ is connected and $f(x) > f(\bar x)$ whenever $\bar x \ne x \in U$.  But the image $f(U)$ must then be connected, and hence contained in the interval $(-\infty,f(\bar x)]$.  We have arrived at the contradiction $U = \{\bar x\}$.  Since $\M$ is strongly critical, we know 
$\overline{|\nabla f|} (y_i) \to 0$, contradicting inequality (\ref{big-slope}).  Property (\ref{small-slope}) follows.

We note that the identifiability of $\M$ ensures that any talweg $x: (f(\bar x), r) \to X$ lies in $\mathcal{M}$ for all $r$ near $f(\bar x)$. This correlates with our observation that the KL property is determined just by the behavior of $f$ on $\mathcal{M}$.

\part{Identifiability in Euclidean spaces}\label{part:2}
\section{Four classes of nonsmooth functions} \label{four}
The second part of this work focuses on identifiable sets in Euclidean space.  We focus initially on four objective classes, consisting of functions that are, in increasing order of strength, prox-regular, primal lower nice, strongly amenable, or ``smoothly embedded''.  The first three are broad classes, familiar from earlier literature \cite{VA,poliquin-pln}.  The last class, consisting just of smooth functions, indicator functions of manifolds, and their sums, is much narrower.  Nonetheless, we shall see in Part II that, from the perspective of identifiability and its continuous-time analog, smoothly embedded functions often model the behavior of much more general functions. 

We begin by reviewing standard variational analysis \cite{VA}.  Throughout Part II, {\bf the term ``smooth''  means ${\mathcal C}^{(2)}$-smooth}.   For a function $f \colon \Rn \to (-\infty,+\infty]$ and a point $\bar x \in \mathrm{dom}\, f$, a vector $y \in \R^n$ is a {\em regular subgradient} of $f$ at $\bar x$, written $y \in \hat \partial f(\bar x)$, if $f(x) \ge f(\bar x) + \langle y, x - \bar x \rangle + o(|x - \bar x|)$ as $x \to \bar x$.  If furthermore we can replace the term $o(|x - \bar x|)$ by the term $O(|x - \bar x|^2)$, then we instead call $y$ a {\em proximal subgradient}, written $y \in \partial_P f(\bar x)$. We say that $y$ is a {\em subgradient} of $f$ at $\bar x$, written $y \in \partial f(\bar x)$, if there are sequences $x_r \to_f \, \bar x$ and $y_r \in \hat \partial f(x_r)$ with $y_r \to y$.  As we observed after Proposition \ref{distance-formula}, $\bar x$ is critical if and only if $0 \in \partial f (\bar x)$.  By contrast, it is easy to check that $f$ has slope zero at $\bar x$ if and only if $0 \in \hat\partial f (\bar x)$.  We define the {\em regular normal cone} of a set $\M \subset \R^n$ at a point $\bar x \in \M$ by $\widehat{N}_\M (\bar x) = \hat \partial \delta_{\M} (\bar x)$, and the {\em normal cone} by $N_\M (\bar x) = \partial \delta_{\M} (\bar x)$. When $\M$ is convex, $\mathrm{ri}\, \M$  denotes its relative interior, and
$\mathrm{par}\, \M$ denotes the unique linear subspace that is a translate of its affine span. 
We denote the closed unit ball by $B$.

We next summarize some basic ideas about manifolds:  more details may be found in standard references such as \cite{Lee}.  In parallel with our terminology for functions, throughout Part II, {\bf the term ``manifold'' means a ${\mathcal C}^{(2)}$-smooth manifold}.  Given a manifold $\M \subset \Rn$ and a point $\bar x \in \M$, there exists a smooth map \mbox{$F \colon \Rn \to \R^m$} with $F(\bar x) = 0$ and surjective Jacobian $\nabla F(\bar x)$ such that 
\bmye \label{zero}
x \in \M \quad \Leftrightarrow \quad F(x) = 0 \qquad \mbox{for all $x \in \Rn$ near $\bar x$}.
\emye
Consequently \cite[Example~6.8]{VA} both the normal cone $N_\M(\bar x)$ and the regular normal cone $\widehat N_\M(\bar x)$ coincide with the classical normal space to $\M$ at $\bar x$.  The nearest-point projection $\mbox{Proj}_{N_\M(x)}(y)$ is well-defined and depends continuously on points $x,y$ near $\bar x$ with $x \in \M$.  A function 
\mbox{$h \colon \M \to \R$} is {\em smooth} around 
$\bar x$ if there exists a smooth function $\tilde h \colon \Rn \to \R$ such that $h(x) = \tilde h(x)$ for all points $x \in \M$ near $\bar x$. We then call $\tilde h$ a {\em smooth extension} of $h$ at $\bar x$, and the {\em Riemannian gradient} of $h$ at $\bar x$, denoted 
$\nabla_\M h(\bar x)$, is the orthogonal projection of the gradient $\nabla \tilde h(\bar x)$ onto the tangent space $T_\M(\bar x)$.  The projected vector is in fact independent of the choice of extension $\tilde h$.  For example, the set $\M \subset \R^2$ defined by the equation $x_2 = x_1^2$ is a manifold, and the restriction $f|_\M$ of the function $f$ defined by equation $(\ref{simple1})$ is smooth because it agrees with the smooth function $\tilde f(x) = 5x_2 - 4x_1^2$ on $\M$.  At the point $0 \in \R^2$,
the Riemannian gradient $\nabla_\M (f|_\M)$, which we abbreviate to $\nabla_\M f$, is zero.

We consider four objective classes of interest, beginning with prox-regularity \cite{VA}.

\begin{defn} \label{def:prox}
{\rm
A closed function $f \colon \Rn \to (-\infty,+\infty]$ is {\em prox-regular} at a point $\bar x$ for a subgradient $\bar y \in \partial f(\bar x)$ if there exists a constant $\rho \ge 0$ such that all points $x, x'$ near $\bar x$ with $f(x)$ near $f(\bar x)$ and subgradients $y \in \partial f(x)$ near 
$\bar y$ satisfy the inequality $f(x') \ge f(x) + \langle y, x' - x \rangle - \rho |x' - x|^2$.
}
\end{defn}  

\noindent
The second class \cite{poliquin-pln} is clearly more restrictive than prox-regularity.  

\begin{defn} \label{def:pln}
{\rm
A closed function $f \colon \Rn \to (\infty,+\infty]$ is {\em primal lower nice} at a point $\bar x \in \Rn$ if there exist constants $\alpha$ and $\beta$ such that all points 
$x,x'$ near $\bar x$, and subgradients $y \in \partial f(x)$ satisfy
\bmye \label{pln}
f(x') ~\ge~ f(x) ~+~ \ip{y}{x'-x} ~-~ (\alpha + \beta|y|)|x'-x|^2.
\emye
If this property holds for all $\bar x \in \mbox{dom}\, f$, then we simply call $f$ {\em primal lower nice}.
}
\end{defn}

\noindent
The third class, discussed in \cite[Definition 10.23 and Proposition 13.32]{VA}, while versatile, is more restrictive still \cite[Theorem~5.1]{poliquin-pln}.

\begin{defn} \label{amenable}
{\rm
At a point $\bar x \in \Rn$, a function $f$ is {\em strongly amenable}  if it has the local representation $f=g \circ F$ for some closed convex function $g \colon \R^m \to (-\infty,+\infty]$ and smooth map $F \colon \Rn \to \R^m$ such that $F(\bar x)$ lies in $\mbox{dom}\, g$ and the normal cone to $\mbox{cl}(\mbox{dom}\, g)$ there intersects the null space of the adjoint map $\nabla F(\bar x)^*$ trivially.
}
\end{defn}

\noindent
In the scenario of Definition \ref{amenable}, a standard chain rule \cite{VA} implies
\[
\partial f(\bar x) ~=~ \hat\partial f(\bar x) ~=~ \nabla F(\bar x)^* \partial g\big( F(\bar x)\big).
\]

\noindent
The final class can be seen, via property (\ref{zero}), to consist of particularly simple strongly amenable functions.

\begin{defn} \label{embedded}
{\rm
At a point $\bar x \in \Rn$, a function $f$ is {\em smoothly embedded}  if it has the local representation $f = g + \delta_\M$ for some smooth function $g \colon \Rn \to \R$ and some manifold $\M$ containing $\bar x$.
}
\end{defn}

\noindent
In the scenario of Definition \ref{embedded}, a standard sum rule \cite[Corollary 10.9]{VA} implies
\bmye \label{restriction}
\partial f(\bar x) ~=~ \hat\partial f(\bar x) ~=~ \nabla g(\bar x) + N_\M(\bar x).
\emye

%
%

The function (\ref{simple1}) is strongly amenable at zero, because we can represent it using the smooth map $F(x) = (x_2-x_1^2,x_1)$ and the convex function $g(u,v) = 5|u|+v$.  
{\em Weakly convex} functions --- those of the form $g - \rho|\cdot|^2$, for proper closed convex functions $g$ and constant $\rho$ (see for example \cite{dd-weakly}) --- are easily seen to be strongly amenable.  For locally Lipschitz functions, prox-regularity, the primal lower nice property, and strong amenability  are equivalent, each coinciding with the ``lower ${\mathcal C}^2$ property'' \cite{VA}.  Indeed, a locally Lipschitz $f$ has each of these properties at a point if and only if it equals a weakly convex function locally \cite[Exercise 10.36, Proposition 13.33]{VA}. 

In general, however, prox-regular functions may not be primal lower nice, even if they are continuous.  For example,  the function $f(x) = \sqrt{|x|}$ is prox-regular at zero, but not primal lower nice there.  To see this, assume that inequality (\ref{pln}) holds.  For any small
$t>0$, we can set $x = t^2$, $x' = 4t^2$, $y = \frac{1}{2t}$, 
to deduce
\[
2t ~\ge~ t + \frac{1}{2t} 3t^2 - \Big(\alpha + \beta\frac{1}{2t}\Big)(3t^2)^2.
\]
This gives a contradiction for $t$ sufficiently small.  Also worth noting is that primal lower nice functions may not be weakly convex:  a smoothly embedded function as described in Definition \ref{embedded} is weakly convex if and only if the set $\M$ is locally affine.

We prove two simple tools for later use.  First, the primal lower nice property is preserved by addition of smooth functions.

\begin{prop}[Primal lower nice preservation] \label{preservation} 
For a function
$g \colon \Rn \to \R$ that is smooth around the point $\bar x \in \Rn$, if a closed function $f \colon \Rn \to (-\infty,+\infty]$ is primal lower nice at $\bar x$, then so is the sum $f+g$.
\end{prop}

\pf
Using the notation above, we can assume that inequality (\ref{pln}) holds.  Since $g$ is smooth, its gradient is Lipschitz around $\bar x$.  (In fact this property suffices for our proof.) Hence, for some constants $\gamma, \lambda>0$, all points $x$ and $x'$ near $\bar x$ satisfy
\bmye \label{Lsmooth}
g(x') ~\ge~ g(x) ~+~ \ip{\nabla g(x)}{x'-x} ~-~ \gamma|x'-x|^2
\emye
and $|\nabla g(x)| \le \lambda$.

Now consider the sum $h=f+g$.  For any points $x$ and $x'$ near $\bar x$ and subgradient 
$w \in \partial h(x)$, the sum rule \cite[Corollary 10.9]{VA} ensures the existence of a subgradient 
$y \in \partial f(x)$ satisfying $w = y + \nabla g(x)$.  Adding the inequalities (\ref{pln}) and 
(\ref{Lsmooth}), we deduce
\[
h(x') 
~\ge~ h(x) ~+~ \ip{w}{x'-x} ~-~ \big(\gamma +  (\alpha + \beta|y|)\big) |x'-x|^2.
\]
Since $|y| \le \lambda + |w|$, the primal lower nice property for $h$ follows.
\finpf

Secondly, we show that when considering primal lower nice functions, we lose no generality in assuming that they are bounded below.

\begin{lem}[Localization] \label{restrict}
If a closed function $f \colon \Rn \to (-\infty,+\infty]$ is primal lower nice at a point $\bar x \in \Rn$, then there is another proper closed function that is primal lower nice and bounded below, and that agrees identically with $f$ near $\bar x$.
\end{lem}

\pf
For simplicity, suppose $\bar x = 0$.  By definition, $f$ is primal lower nice throughout the ball 
$2\delta B$ for some $\delta > 0$.  Since $f $ is closed, we can shrink $\delta > 0$ if necessary to ensure $f$ is bounded below on $2\delta B$.  The function $g \colon \Rn \to \oR$ defined by
\[
g(x) ~=~
\left\{
\begin{array}{ll}
0 & (|x| \le \delta) \\
\frac{(|x|^2 - \delta^2)^3}{4\delta^2 - |x|^2} & (\delta < |x| < 2\delta) \\
+\infty & (|x| \ge 2\delta)
\end{array}
\right.
\]
is smooth on the interior of $2\delta B$, so, by Proposition \ref{preservation}, the function $f+g$ is primal lower nice throughout its domain.  But $f+g$ is also bounded below and agrees identically with $f$ on the ball $\delta B$, as required.
\finpf

A particularly simple class of strongly amenable functions $f$ arise when the convex function $g$ in Definition \ref{amenable} is given by $g(y) = \max_i y_i$.  The resulting functions $f$ illustrate well the theory of identifiability:  we discuss them in the next section.

\section{Max functions} \label{max}
We consider the important basic example of {\em max functions\/}:  those functions of the form $x \mapsto \max_i f_i(x)$ for $x \in \Rn$, where the index $i$ ranges over some finite index set and each function $f_i$ is continuously differentiable.  Identifiable sets for max functions were studied in \cite[Corollary 4.8]{Lewis-active}, but under nondegeneracy conditions that often fail in practice \cite{mirror}.  By contrast, \cite[Corollary 4.11]{ident} dispenses with nondegeneracy, but omits important details in the result and proof.  Here we present a short and direct argument.  We start with the smooth case.

\begin{prop} \label{differentiable}
If a function $f$ on $\Rn$ is continuously differentiable around a point 
$\bar x$, with $\nabla f(\bar x) = 0$, then the intersection of any set containing $\bar x$ with any sufficiently small neighborhood of $\bar x$ is strongly critical at $\bar x$.
\end{prop}

\pf
As noted, $|\nabla f|(x) = |\nabla f(x)|$ at all points $x$ near $\bar x$.  The result follows.
\finpf

The next result uses Proposition \ref{distance-formula} and standard variational analysis: \cite[Theorem 10.31]{VA} reveals the property $G(x) = \partial f(x)$.  A direct proof is also short.

\begin{prop}[Slope for max functions] \label{setting}
Consider a nonempty finite index set $I$ and continuously differentiable functions 
$f_i \colon \Rn \to \R$ for $i \in I$.  For each point $x \in \Rn$, define the value $f(x) = \max_i f_i(x)$, and sets
\[
I(x) ~=~ \{i : f_i(x) = f(x)\} \qquad \mbox{and} \qquad
G(x) ~=~ \mbox{\rm conv} \{  \nabla f_i(x) : i \in I(x)\}.
\]
Then the slope and limiting slope of $f$ at $x$ both equal $\mbox{\rm dist}\big(0,G(x)\big)$.
\end{prop}

\pf
Elementary calculus shows
\[
|\nabla f|(x) ~=~ \max_{|y|=1} \min_{i \in I(x)} \ip{\nabla f_i(x)}{y}.
\]
A simple limiting argument then shows that the slope and limiting slope agree everywhere, and the result then follows from the separating hyperplane theorem.
\finpf

For any index set $I$, the {\em support set} of a vector $\lambda \in \R^I$ is the set
$\{i \in I : \lambda_i \ne 0 \}$.  

\begin{thm}[Identifiable sets for max functions] \label{max-thm}
With the assumptions of Proposition \ref{setting}, consider a critical point 
$\bar x$ for the max function $f$.  Let $\bar I = I(\bar x)$, and denote by $\Omega$ the finite set consisting of support sets of vectors $\lambda \in \R_+^{\bar I}$ satisfying 
\bmye \label{Lambda}
\sum_{i \in \bar I} \lambda_i = 1, \qquad \sum_{i \in \bar I} \lambda_i \nabla f_i(\bar x) = 0.
\emye
Then the union of the sets $\M_J = \{x : I(x) \supset J\}$ as the support set $J$ ranges over $\Omega$ 
is strongly critical and identifiable at $\bar x$.
\end{thm}

\pf
We first prove identifiability.  Consider a sequence $x_r \to \bar x$ satisfying 
$|\nabla f|(x_r) \to 0$ and with $x_r \not\in \M_J$ for $r=1,2,\ldots$ and all $J \in \Omega$.  For large $r$, continuity implies $I(x_r) \subset \bar I$, and  Proposition \ref{setting} then implies the existence of a vector $\lambda^r \in \R^{\bar I}_+$ satisfying $\sum_i \lambda^r_i = 1$ and  $\lambda^r_i = 0$ for all $i \not\in I(x_r)$, with $\sum_i \lambda^r_i \nabla f_i(x_r) \to 0$. After taking a subsequence, $\lambda^r$ converges to some vector $\bar\lambda \in \R_+^{\bar I}$ with support set $\bar J \subset \bar I$.  Since $\bar\lambda$ satisfies equation (\ref{Lambda}), by continuous differentiability, $\bar J \in \Omega$.  Consequently $I(x_r) \not\supset \bar J$ for each $r$, and hence there exists an index 
$i(r) \not\in I(x_r)$ with $\bar\lambda_{i(r)} > 0$.  Taking a further subsequence, some index $i \in \bar I$
satisfies $\bar\lambda_i > 0$ and $i \not\in I(x_r)$ for all $r$.  This gives the desired contradiction: 
$0 = \lambda^r_i \to \bar\lambda_i > 0$.

Turning to strong criticality, by Proposition \ref{finite}, it suffices to prove, for each support set 
$J \in \Omega$, that the set $\M_J$ is strongly critical.  Some corresponding solution $\lambda \in \R_+^{\bar I}$ of equation (\ref{Lambda}) has support set $J$.  The continuously differentiable function $g = \sum_i \lambda_i f_i$ satisfies the inequality $g \le f$, with equality throughout $\M_J$.  Strong criticality follows by Propositions \ref{strongly} and \ref{differentiable}.
\finpf

To illustrate, consider the max function on $\R^2$  with the two ingredients defined by $(u,v) \mapsto v^2 \pm u$ at the critical point $(0,0)$.   Theorem \ref{max-thm} then recovers Example~\ref{basic}.  Although in this example the strongly critical identifiable set is a manifold (defined by $u=0$), that might not be true more generally.  For example, the function on $\R$ with ingredients the zero function and the function $x \mapsto -x$, at the critical point $0$, has $\R_+$ as a strongly critical identifiable set.  Furthermore, in general, the set $\Omega$ in Theorem \ref{max-thm} may contain several distinct supports.  Consider, for example, the function on $\R$ with ingredients the zero function and the function $x \mapsto x^3$, at the critical point $0$, in which case $\R$ is a strongly critical identifiable set.  When the gradients $\nabla f_i(\bar x)$ for $i \in \bar I$ are affinely independent, however, equation (\ref{Lambda}) has a unique solution $\bar\lambda$.   If we fix some index $i$ in the index set $I_> = \{i \in \bar I : \bar\lambda_i > 0 \}$, then $M$ is defined around $\bar x$ by the constraints
\begin{eqnarray*}
f_j(x) - f_i(x) &=& 0 \qquad \mbox{for}~ j \in I_> \setminus \{i\} \\
f_j(x) - f_i(x) &\le& 0 \qquad \mbox{for}~ j \in \bar I \setminus I_>.
\end{eqnarray*}
The functions on the left-hand side are continuously differentiable, with linearly independent gradients, by the affine independence assumption.  The inverse function theorem therefore guarantees that, locally, the set $\mathcal{\M}$ is diffeomorphic to a polyhedral cone.  Consequently, around $\bar x$, the strongly critical identifiable set $\mathcal{M}$ is normally embedded, so Theorem \ref{linear-proj} applies.

\section{Identifiability and partial smoothness}\label{sec:ps}
In Euclidean space, identifiability, an idea inspired by explorations in \cite{Flam92,Burke-More88,Burke90,Wright}, was introduced in \cite[Definition 3.10]{ident}:  by Corollary~\ref{original}, that original definition coincides with ours. A primary aim of \cite{ident} was to relate identifiability and the idea of {\em partial smoothness} introduced in \cite{Lewis-active}.  We now revisit that relationship using an fresh proof scheme, replacing the intricate proximal epigraphical analysis suggested in \cite{ident} by a direct and natural argument based on Theorem \ref{linear} (Linear growth).
The next result, following a pattern observed in\cite[Section 10]{ident2}, is the key tool.

\begin{prop}[Identifiable sets and subgradients] \label{lem:regular2}
Suppose that the closed function $f \colon \Rn \to (-\infty,+\infty]$ has slope zero at the point $\bar x$.  Consider a closed, normally embedded set $\M \subset \Rn$ that is strongly critical and identifiable at $\bar x$, and has continuous restriction $f|_\M$.  Then for all points $x$ near $\bar x$ with value $f(x)$ near $f(\bar x)$ and small vectors $y \in \Rn$, the following two properties are equivalent:
\begin{itemize}
\item[{\rm (a)}]
$y \in \partial f(x)$
\item[{\rm (b)}]
$x \in \M$ and $y \in \partial (f+\delta_\M)(x)$.
\end{itemize}
The same equivalence holds for both regular and proximal subgradients.  Furthermore, prox-regularity at $x$ for $y$ holds for $f$ if and only if it holds for $f+\delta_\M$.
\end{prop}

\pf
We first prove the regular subdifferential result.  In that case, property (a) implies $x \in \M$, by identifiability, and then property (b) follows from the inclusion 
$\hat\partial f(x) \subset \hat\partial (f+\delta_\M)(x)$. Hence, in the regular case, we just need to prove (b) implies (a).  By Theorem \ref{linear-proj}, there exists $\epsilon > 0$ such that all points 
$x \in B^f_{2\epsilon}(\bar x)$ and nearest points $v \in \mbox{Proj}_\M(x)$ satisfy 
$f(v) + \epsilon|v-x| \le f(x)$.  To complete the proof, we claim 
\[
x \in \M \cap B^f_{\epsilon}(\bar x), \quad y \in \hat\partial(f+\delta_\M)(x), \quad \mbox{and} 
\quad |y|<\epsilon \quad \Rightarrow \quad y \in \hat\partial f(x).
\]
Assume $x$ and $y$ satisfy the left-hand side properties.  Then points $v \in \M$ satisfy
\[
f(v) ~\ge~ f(x) + \ip{y}{v-x} + o(|v-x|) \qquad \mbox{as}~ v \to x.
\]
If $y \not\in \hat \partial f(x)$, then there exists $\delta > 0$ and a sequence $z_s \to x$ satisfying 
\[
f(z_s) - f(x) - \ip{y}{z_s-x} ~<~ -\delta |z_s-x|.
\] 
For all large $s$, since $f$ is closed, $f(z_s) > f(x) -\epsilon$, and we also have 
$f(z_s) < f(x)+ \epsilon$, so in fact 
$z_s \in B_\epsilon^f(x) \subset B_{2\epsilon}^f(\bar x)$.  Choosing  $v_s \in \mbox{Proj}_\M(z_s)$, we have
\begin{eqnarray*}
f(z_s) & \ge & f(v_s) + \epsilon|v_s-z_s| \\
& \ge & f(x) + \ip{y}{v_s-x} + \epsilon|v_s-z_s| + o(|v_s-x|) \\
& = &
f(x) + \ip{y}{z_s-x} + \ip{y}{v_s-z_s} + \epsilon|v_s-z_s| + o(|v_s-x|) \\
& \ge & f(x) + \ip{y}{z_s-x} + o(|v_s-x|),
\end{eqnarray*}
so $-\delta |z_s- x| > o(|v_s-x|)$, and hence in fact $|z_s- x| = o(|v_s-x|)$, contradicting
\[
|v_s-x| \le |v_s-z_s| + |z_s-x| \le 2|z_s-x|.
\]
This completes the proof in the regular case.  The analogous proof for proximal subgradients is almost identical, and the result for subgradients follows from the regular case by a simple limiting argument.  Finally, prox-regularity of $f$ at $x$ for a proximal subgradient $y$ depends only on the behavior of $f$-attentive localizations of the mapping $\partial f$ around $(x,y)$ \cite[Theorem 13.36]{VA}, so the final claim follows.
\finpf

\begin{defn} \label{def:ps}
{\rm
Given a closed function $f \colon \Rn \to \overline{\R}$, consider a point $\bar x \in \mbox{dom}\, f$ and a manifold $\M$ containing $\bar x$. Then $f$ is {\em partly smooth} at $\bar x$ for a subgradient 
$\bar y \in \partial f(\bar x)$ relative to $\M$ if it satisfies the following properties:
\begin{quote}
\begin{description}
\item[Prox-regularity:]  the function $f$ is prox-regular at $\bar x$ for $\bar y$.
\item[Restricted smoothness:]  the restriction $f|_{\M}$ is smooth around $\bar x$.
\item[Sharpness:] the subspace $\mathrm{par}\, \hat\partial f(\bar x)$ is just the normal space
$N_{\M} (\bar x)$.
\item[Inner semicontinuity:] for all $y \in \partial f(\bar x)$ near $\bar y$ and
sequences \\ $x_r \to \bar x$ in $\M$, there exist $y_r \in \partial f(x_r)$ converging to $y$.
\end{description}
\end{quote}
}
\end{defn}

\noindent
This definition, from \cite{ident2}, slightly refined the original \cite{Lewis-active}.
We note an easy example.

\begin{prop} \label{smoothly-embedded}
Any smoothly embedded function is partly smooth at every point in its domain $\M$ relative to $\M$.
\end{prop}

\pf
As we have observed, any smoothly embedded function $f=g+\delta_\M$, for a smooth function $g$ and a manifold $\M$, must be prox-regular.  The restricted smoothness condition is immediate.  The sharpness condition follows from formula~(\ref{restriction}).  Finally, given any subgradient 
$y \in \partial f(\bar x)$, define a vector $\tilde y = y - \nabla g(\bar x) \in N_\M(\bar x)$.  Given any sequence $x_r \to \bar x$ in $\M$, we have
\[
\nabla g(x_r) + \mbox{Proj}_{N_\M(x_r)}(\tilde y)
~\to~ 
\nabla g(\bar x) + \mbox{Proj}_{N_\M(\bar x)}(\tilde y) ~=~ \nabla g(\bar x) + \tilde y ~=~ y.
\]
The left-hand side vector lies in $\partial f(x_r)$, so inner semicontinuity follows.
\finpf

We focus on partial smoothness for the subgradient $\bar y = 0$.  By Theorem \ref{distance-formula}, the following result follows immediately from the inner semicontinuity condition.

\begin{prop}
If a closed function $f \colon \Rn \to (-\infty,+\infty]$ is partly smooth for zero at a point $\bar x$ relative to a manifold $\M$, then $\M$ is strongly critical at $\bar x$.
\end{prop}

Although the properties required for partial smoothness appear at first to comprise a complex list, remarkably, in semi-algebraic optimization, partial smoothness holds generically \cite{gen_exist}, along with the {\em nondegeneracy condition} $0 \in \mbox{ri}\,\hat\partial f(\bar x)$.
More precisely, for all vectors $y$ in a dense open semi-algebraic subset of $\Rn$, the perturbed function $f + \ip{y}{\cdot}$ is partly smooth at all critical points $x$, each of which furthermore must be nondegenerate.  Consequently, the perturbed function has an identifiable manifold at each critical point, for zero.  Thus, for generic semi-algebraic optimization problems, we expect to see identifiable manifolds around the solutions.  Partial smoothness and nondegeneracy also imply an important constant-rank property for the subdifferential operator $\partial f$ \cite[Theorem~6.5]{Lewis-Liang-Tian}.  

Most striking from our current perspective is the following result \cite[Proposition~10.12]{ident2}:  partial smoothness and nondegeneracy together are equivalent to a seemingly simpler identifiability property with respect to a manifold.  The argument is now a straightforward application of Proposition \ref{lem:regular2}.

  \begin{thm}[Identifiability and partial smoothness] \label{thm:ps-ident}~~
Consider a closed function $f \colon \Rn \to (-\infty,+\infty]$ with slope zero at a point $\bar x$ lying in a manifold $\M$.  Suppose that $f|_{\M}$ is smooth around $\bar x$. Then $\M$ is identifiable at 
$\bar x$ if and only if $f$ is partly smooth at $\bar x$ for zero relative to $\M$ with 
$0 \in \mathrm{ri}\, \hat \partial f(\bar x)$.
\end{thm}

\pf
Suppose that the manifold $\M$ is identifiable at the point $\bar x$.  By assumption, the restriction $f|_\M$ has a smooth extension $\tilde f$ at $\bar x$.  Like the Euclidean case, it is easy to verify that the slope $|\nabla (f|_\M)|(x)$ is just given by $|\nabla (f|_\M)(x)|$, at all points $x$ near $\bar x$, and hence is continuous at $\bar x$. Applying Theorem \ref{linear-proj}, we know $\M$ is strongly critical at $\bar x$.  Consider the function 
$h = f+\delta_\M = \tilde f + \delta_\M$.  We have 
\[
0 ~\in~ \hat\partial f(\bar x) ~\subset~ \hat\partial h(\bar x) ~=~ 
\nabla\tilde f(\bar x) + N_\M(\bar x),  
\]
so $\nabla\tilde f(\bar x) \in N_\M(\bar x)$ and $\hat\partial h(\bar x) = N_\M(\bar x)$.
Using Proposition \ref{lem:regular2}, the convex sets $\hat\partial h(\bar x)$ and $\hat \partial f(\bar x)$ are identical near zero so the nondegeneracy condition $0 \in \mbox{ri}\,\hat\partial f(\bar x)$ and the sharpness condition follow, and since $h$ is smoothly embedded, it is prox-regular at $\bar x$ for $0$, whence so is $f$.  Restricted smoothness holds by assumption.

Inner semicontinuity also holds by Proposition \ref{lem:regular2}.  To see this claim, note that any small subgradient $y \in \partial f(\bar x)$ also lies in 
$\partial h(\bar x)$.  Given any sequence $x_r \to \bar x$ in $\M$, since the smoothly embedded function $h$ is partly smooth by Proposition~\ref{smoothly-embedded}, there exist subgradients 
$y_r \in \partial h(x_r)$ approaching $y$.  For all large $r$ we have \mbox{$y_r \in \partial f(x_r)$}, so the claim follows.

Our proof in the converse direction is more standard, following \cite{Hare} but shorter.  If the result fails, then by Corollary \ref{original} there exists a sequence of points $x_r \to_f \bar x$ in $\M^c$ with subgradients $y_r \in \partial f (x_r) \to 0$.  For large $r$, the point $x_r$ has a unique nearest point $x'_r$ in $\M$.  Normalizing the vector 
$x_r-x'_r$ gives a unit vector $w_r \in N_\M(x'_r)$, which, taking a subsequence, we can suppose converges to a unit vector $w \in N_\M(\bar x)$.  The sharpness and nondegeneracy conditions imply the existence of a constant $\epsilon > 0$ such that $\epsilon w \in \partial f(\bar x)$, and the continuity condition then implies the existence of a sequence of subgradients $y_r' \in \partial f(x_r') \to \epsilon w$. Since $f$ is prox-regular at $\bar x$ for $0$, there exists a constant $C > 0$ such that, for all large $r$ we have
\begin{eqnarray*}
f(x_r') 
& \ge & f(x_r) + \langle y_r, x_r' - x_r \rangle - C |x_r' - x_r|^2 \quad \mbox{and} \\
f(x_r) 
& \ge & f(x_r') + \langle y_r', x_r - x_r' \rangle - C |x_r' - x_r|^2,
\end{eqnarray*}
and hence $\ip{y'_r-y_r}{w_r} \le 2C|x_r-x'_r|$.  But in this inequality, the left-hand side converges to $\epsilon$, contradicting the convergence of the right-hand side to zero.
\finpf

\noindent
The original result, \cite[Proposition 10.12]{ident2}, concerns partial smoothness for an arbitrary subgradient, rather than zero specifically.  This more general result follows immediately by adding a linear perturbation to the objective $f$.

\subsection*{Partly analytic functions}
We end this section by discussing an important special case of Theorem~\ref{thm:ps-ident}.
We present a simple tool that often suffices for verifying the KL property of Section~\ref{KL} and helps find desingularizers.  We describe how the KL property can decouple, becoming, like partial smoothness, a simple and intuitive hybrid of smooth and sharp behavior.  In general, semi-algebraic (or tame) functions always have the KL property, but the proof \cite{Lewis-Clarke} is not straightforward, revealing little about desingularizers.  For {\em generic} semi-algebraic objectives, on the other hand --- specifically for the sum of any specific semi-algebraic function and a generic linear perturbation --- the situation is simpler:  \cite{gen,gen_exist} prove the existence of an identifiable analytic manifold on which the objective is analytic.  We call this situation the {\em partly analytic} case.

We can equip any manifold $\M$ in Euclidean space with two natural metrics:  the metric $d$ induced by the Euclidean distance, and  the Riemannian distance $d_{\mbox{\rm\scriptsize Rie}}$. For any points $x,y \in \M$, $d_{\mbox{\rm\scriptsize Rie}}(x,y)$ is the infimum of the length of paths in $\M$ from $x$ to $y$ \cite{lee-riemannian}.  We can move easily between these metrics using the following technique.

\begin{prop} \label{infinitesimal}
Given a manifold $\M$ in Euclidean space, consider a function $g \colon \M \to (-\infty,+\infty]$, a point $\bar x \in \mbox{\rm dom}\,g$, and a desingularizer $\phi$.  Then, with respect to the Riemannian and induced Euclidean metrics $d_{\mbox{\rm\scriptsize Rie}}$ and $d$, the KL properties for $g$ at $\bar x$ with desingularizer $\phi$ are equivalent.
\end{prop}

\pf
This is a consequence of the following fact, proved in \cite[Proposition 3.1]{dontchev-lewis}:
\[
\mbox{} \hspace{4cm} \lim_{x,y \to \bar x \atop x \ne y} \frac{d_{\mbox{\rm\scriptsize Rie}}(x,y)}{d(x,y)} ~=~ 1. \hspace{4cm} \Box
\]
In the partly analytic case, the KL property follows from the following result.  The result does not in fact require tameness. For example, it applies to the objective  $x \mapsto f(x) + |x|$, for the function $f$ defined in Example \ref{unidentifiable}.  This objective is not tame, but does satisfy the KL property at zero, since $\{0\}$ is an identifiable set.

\begin{cor}[Partly analytic functions]
Suppose that $f \colon \Rn \to (-\infty,+\infty]$ is closed and has slope zero at a point $\bar x$.  If some analytic manifold $\M$ is identifiable at $\bar x$ and the restriction $f|_\M$ is analytic, then the KL property for $f$ holds at $\bar x$, with $KL$ exponent equal to that of $f|_\M$ at $\bar x$.
\end{cor}

\pf
Analytic functions on analytic manifolds satisfy the KL property \cite[Section~9]{kurdyka-mostowski-parusinski}.  The conclusion then follows from Theorem \ref{kl-identifiable} and Proposition \ref{infinitesimal}.
\finpf

\section{Subgradient curves and identification} \label{sec:Subgradient}
The property of identifiability fundamentally concerns critical sequences.  We now move on from that discrete-time perspective to study continuous-time trajectories. Consider once again the question of minimizing an objective function $f$ on a metric space.  The central continuous-time analogue of the sequences considered thus far is the family of curves of ``maximal slope'' \cite{giorgi,ambrosio}.  Around some given point $x^*$, we pose the question of the existence of a small set $\M$ such that curves of maximal slope converging to $x^*$ always remain in $\M$ eventually.

Figure \ref{fig:subcurve} illustrates this behavior for the simple nonsmooth nonconvex objective $5 |x_2-x_1^2| + x_1^2$,  around the minimizer $x^* = 0$:  the associated small set 
$\M$  is defined by $x_2=x_1^2$.  The behavior is the precise analogue of the phenomenon we observed in preceding sections, where various iterative procedures for minimizing $f$ identify an associated set $\M$ around $x^*$:  sequences of iterates converging to $x^*$ must eventually lie in $\M$, thereby revealing some solution structure like constraint activity, sparsity, or matrix rank.  For objectives $f$ on general metric spaces, the question of identification in continuous time remains open.  As a first step, therefore, we focus on the Euclidean case, $f \colon \Rn \to (-\infty,+\infty]$.  Following the arguments of \cite{curves},  for interesting functions $f$ the question then typically reduces to the behavior of {\em subgradient curves\/}:  locally absolutely continuous maps 
$x \colon \R_+ \to \mbox{dom}\, f$ satisfying
\[
x'(t) \in -\partial f\big(x(t)\big) \qquad \mbox{for almost all times}~ t>0.
\]
Our main continuous-time result, Theorem \ref{main-1}, demonstrates for the first time identification behavior for subgradient curves.

The significance of subgradient curves from an optimization perspective has been highlighted in 
\cite{attouch-teboulle,attouch-bolte}.  Existing literature covers a wide class of nonsmooth and nonconvex functions $f$ for which they must converge.  In particular, the existence and uniqueness of subgradient curves, given any initial point $x(0)$ in $\mbox{dom}\, f$, was shown for convex functions $f$ in a well-known 1973 monograph of Br{\'e}zis \cite{brezis}.  In fact this result holds more generally, under the assumption that $f$ is primal lower nice, as shown in the following theorem from
\cite{marcellin-thibault}.

\begin{thm}[Existence and uniqueness of subgradient curves] 
\label{thibault} \hfill \mbox{}
If the closed function $f:\Rn \to (-\infty,+\infty]$ is primal lower nice and bounded below, then there exists a unique subgradient curve $x(\cdot)$ corresponding to any initial point $x(0) \in \mbox{\rm dom}\, f$.  The curve furthermore satisfies
\[
\int_0^{\infty}|x'(t)|^2\,dt ~<~ +\infty.
\]
\end{thm}

For closed convex functions with minimizers, following \cite{brezis}, Bruck proved that subgradient curves converge to minimizers \cite{bruck}. Subgradient curves are known to converge more generally, assuming some version of the KL property, holding in particular for closed semi-algebraic (or subanalytic) functions on bounded domains \cite{loja,Lewis-Clarke}. Such convergence results hold for all closed convex functions, ``lower $\mathcal{C}^2$" functions, and weakly convex functions that are bounded below \cite{tailwag, loja}. We note that these results generalize:  for a primal lower nice function $f$ that is bounded below, if the KL property holds throughout some subgradient curve, then a simple argument using \cite[Theorem 3.2]{marcellin-thibault} shows that the curve must have finite length and hence converge.  We deduce, for example, that all subgradient curves for the function (\ref{simple1}) converge to the minimizer at zero.  To summarize, we can reasonably focus our current study on subgradient curves that converge.

The following is our main continuous-time result.

\begin{thm}[Identification for subgradient curves] \label{main-1} 
Consider a closed function $f \colon \Rn \to (-\infty,+\infty]$ and a subgradient curve $x(\cdot)$ converging to a point $\bar{x}$ at which $f$ is primal lower nice and has a closed identifiable set $\M$. Then, after a finite time, $x(\cdot)$ lies on $\M$.  If furthermore $\M$ is normally embedded and strongly critical at $\bar{x}$, and the restriction $f|_\M$ is continuous, then $x(\cdot)$ becomes a subgradient curve for the function $f + \delta_\M$.
\end{thm}

\noindent
To derive continuous-time identification results of this general type, we must prove that the speed of any subgradient curve always converges (essentially) to zero.  In the convex case, this essential convergence property is well known.  Indeed, the convex special case of Theorem \ref{main-1} is part of the identification folklore, although we do not know a reference.  In the next section, we prove essential convergence in the primal lower nice case, and hence deduce Theorem \ref{main-1}.

In Theorem \ref{main-1}, if the identifiable set $\M$ is a manifold, then we deduce that the dynamic of the subgradient curves is eventually governed by a smooth system. Before stating the general result, we consider a special case.

\begin{prop}[Smooth restrictions] \label{restrictions}
Suppose that the set $\M \subset \Rn$ is a manifold around the point $\bar x$, and that the function $f \colon \Rn \to (-\infty,+\infty]$ is smooth around $\bar x$.  Then, near $\bar x$, subgradient curves $x(\cdot)$ for the function $f + \delta_\M$ are just smooth curves in $\M$ solving the classical differential equation
\bmye \label{Riemman}
x'(t) = -\nabla_\M f \big(x(t)\big). 
\emye
\end{prop}

\pf
By definition, at all times $t$ for which $x(t)$ is near $\bar x$, we must also have $x(t) \in \M$, and furthermore, by the sum rule for subdifferentials \cite{VA},
\[
x'(t) ~\in~ -\nabla f\big(x(t)\big) + N_\M\big(x(t)\big) \qquad \mbox{almost surely}.
\]
At such times $t$ we must have $x'(t) \in T_\M\big(x(t)\big)$ since $x (s) \in \M$ for all time $s$ near $t$. Therefore we deduce equation (\ref{Riemman}).  The smoothness of $x(\cdot)$ now follows from classical smooth initial value theory \cite{diff-mani}.
\finpf

\noindent
The general case asserts the same eventual behavior for a broader class of functions. 

\begin{thm}[Subgradient curves on manifolds] \label{main}
Consider a closed function $f \colon \Rn \to (-\infty,+\infty]$ and a subgradient curve $x(\cdot)$ converging to a point $\bar{x}$ at which $f$ is primal lower nice and has an identifiable manifold $\M$. If $f|_{\mathcal{M}}$ is smooth around $\bar{x}$, then, after a finite time, $x(\cdot)$ becomes a smooth curve on $\M$, satisfying 
\[
x(t) \in \M \qquad \mbox{and} \qquad x'(t) = -\nabla_\M f \big(x(t)\big).
\]
\end{thm}

\pf
Manifolds are locally closed. By Theorem~\ref{main-1}, the curve $x(\cdot)$ eventually lies in $\M$ and becomes a subgradient curve of the function $f + \delta_\M$. Let $\tilde f$ be a smooth extension of $f$. Then $f + \delta_\M = \tilde f + \delta_\M$. The result follows by Proposition~\ref{restrictions}.
\finpf
%
%

\section{Trajectories for primal lower nice functions}\label{sec:pln}
Theorem~\ref{thibault} listed several important properties of subgradient curves for primal lower nice functions.  For our main result, Theorem \ref{main-1}, we rely crucially on one more property.  We prove that velocity vectors $x'(\cdot)$ for bounded subgradient curves $x(\cdot)$ {\em essentially converge to zero\/}:  in other words, for all $\epsilon > 0$, there exists a time $T$ such that for almost all times $t \ge T$ we have $|x'(t)| < \epsilon$.

\begin{thm} \label{convergence}
If a closed function $f \colon \Rn \to (-\infty,+\infty]$ is primal lower nice and bounded below, then the velocity of any bounded subgradient curve for $f$ essentially converges to zero.
\end{thm}

\pf
Consider a bounded subgradient curve $x(\cdot)$.
Suppose $x([0, \infty)) \subset r B$. A standard compactness argument guarantees constants $\alpha, \beta$ and $\delta > 0$ such that the primal lower nice inequality (\ref{pln}) holds for all points $x,x' \in r B$ with $|x - x'| < 2 \delta$. 

By assumption, there exists a full-measure subset $\Omega$ of the interval $[0,+\infty)$ such that 
$x'(t) \in -\partial f\big(x(t)\big)$ for all $t \in \Omega$.
Consider the measure $\rho$ on $\R_+$ with density $|x'(\cdot)|$.
By \cite[Lemma 2.1]{marcellin-thibault}, the primal lower nice property ensures the existence of constants $\mu,\nu > 0$ such that
\[
|x'(s)| ~\le~ |x'(t)| \exp(\mu(s-t) + \nu \rho[t,s])
\]
for all times $t<s$ in $\Omega$ satisfying $x([t, s]) \subset x(t) + \delta B$. 

If the result fails, there is a constant $\epsilon > 0$ and a sequence of times 
$s_j \to + \infty$ in the set $\Omega$ satisfying $|x'(s_j)| > \epsilon$ for all $j$.
Taking a subsequence, we can suppose, for each $j$, the inequality $s_{j+1} - s_j > 1$, from which we deduce either
\bmye \label{pln-ineq1}
|x'(s_j)| ~\le~ |x'(t)| \exp(\mu + \nu \rho[t,{s_j}] ) \qquad \text{for all}~ t \in [s_j - 1,s_j] \cap \Omega,
\emye
or there exists $t \in [s_j - 1,s_j]$ such that
\bmye \label{pln-ineq2}
|x(t) - x(s_j)| ~>~ \delta.
\emye
Assume (\ref{pln-ineq1}) holds.  Then $|x'(t)| \exp(\nu \rho[t,s_j] ) > \epsilon e^{-\mu}$ for all
$t \in [s_j - 1,s_j] \cap \Omega$.
Hence either some time $t$ satisfies $\nu \rho[t,s_j] \ge 1$,
in which case $\rho[s_j-1,s_j]  \ge  \frac{1}{\nu}$,
or there is no such time $t$, in which case all $t \in [s_j - 1,s_j] \cap \Omega$ satisfy the inequality
$|x'(t)| > \epsilon e^{-\mu-1}$, implying $\rho[s_j-1,s_j] \ge  \epsilon e^{-\mu-1}$.
On the other hand, if (\ref{pln-ineq2}) holds, we easily deduce
\[
\rho[s_j-1,s_j]  ~ \ge ~ \rho[t,s_j] ~ \ge ~ |x(t) - x(s_j)| ~ > ~ \delta.
\]
We have thus shown the existence of a constant $\epsilon > 0$ such that $\rho[s_{j-1},s_j] \ge  \epsilon$ for each $j$,
and hence, by H\"older's inequality,
\[
\int_{s_j-1}^{s_j} |x'(\tau)|^2 d\tau ~ \ge ~ \epsilon^2.
\]
But this contradicts the conclusion of Theorem~\ref{thibault}.
\finpf

Theorems \ref{thibault} and \ref{convergence} concern functions that are bounded below.  However, our analysis here is local, so Lemma \ref{restrict} shows that the boundedness assumption involves no loss of generality.  The following result is a local version of Theorem \ref{convergence}.

\begin{thm}[Essential convergence] \label{convergence-local} \mbox{} \\
If a closed function $f \colon \Rn \to (-\infty,+\infty]$ is primal lower nice at a point $\bar x \in \Rn$, and a subgradient curve $x(\cdot)$ for $f$ converges to $\bar x$, then $\bar x$ is critical, the function value $f\big(x(\cdot)\big)$ converges to $f(\bar x)$, and the velocity $x'(\cdot)$ essentially converges to zero.
\end{thm}

\pf
The first two conclusions essentially appear in \cite{marcellin-thibault}, and the third appears in Theorem \ref{convergence}, all under the additional assumption that $f$ is bounded below and primal lower nice throughout its domain.  More generally, we can apply Proposition~\ref{restrict} to replace $f$ by a function that is primal lower nice and bounded below, and that is unchanged near the point $\bar x$.  The original subgradient curve converges to $\bar x$, so it coincides eventually with a subgradient curve for the new function.  
\finpf

The proof of our main continuous-time result is now straightforward.
\bigskip

\noindent
{\bf Proof of Theorem \ref{main-1}}.
By Theorem~\ref{convergence-local}, the function value $f\big(x(\cdot)\big)$ converges to $f(\bar x)$, and there exists a full-measure set $\Omega \subset \R_+$ such that 
$x'(t) \in -\partial f\big(x(t)\big)$ for all $t \in \Omega$
and $x'(t) \to 0$ as $t \to +\infty$ in $\Omega$.
Hence by Corollary \ref{original}, we know $x(t) \in \M$ for all large $t \in \Omega$. Note that since $f$ is primal lower nice at the critical point $\bar x$, it has slope zero there. We then deduce by Proposition \ref{lem:regular2} that $x'(t) \in -\partial (f+\delta_{\M}) \big(x(t)\big)$ for all large 
$t \in \Omega$.
Since the set $\M$ is closed, the subgradient curve $x(\cdot)$, being continuous, eventually lies in $\M$.
\finpf

An alternative approach to Theorem \ref{main} bypasses the essential convergence result.  We proceed as follows. By the square integrability of the velocity $x'$, there exists a sequence of times $t_k \rightarrow \infty$ for which $x(t_k) \in - \partial f(x(t_k))$ for all $k$ and $x'(t_k) \rightarrow 0$. We deduce $x(t_k) \in \mathcal{M}$ for all large $k$, and hence, by \cite[Proposition~10.12]{ident2}, 
$\nabla_{\mathcal{M}} f(y) \in \partial f(y)$ for all points $y \in \mathcal{M}$ near $\bar{x}$. Consequently, for sufficiently large $k$, by the uniqueness of the subgradient curve starting at the point $x(t_k)$, we know that the trajectory $x(t_k + \cdot)$ coincides with the solution 
$\tilde{x}: \R_+ \rightarrow \mathcal{M}$ of the initial value problem
$\tilde{x}'(t) = - \nabla_{\mathcal{M}} f\big( \tilde{x}(t) \big)$ with  $\tilde{x} (0) = x(t_k)$.

\bigskip

\noindent
\textbf{Acknowledgement} The authors thank D.~Drusvyatskiy for helpful suggestions, in particular leading to the alternative argument for Theorem \ref{main}, and two anonymous referees, whose critiques and questions greatly improved the manuscript.


\def\cprime{$'$} \def\cprime{$'$}

\end{document}